\documentclass[12pt]{article}                
\usepackage{graphicx}
\usepackage{psfrag}
\usepackage{amsmath, amssymb}
\usepackage{mathptmx}
\usepackage{ntheorem}

\makeatletter
\newsavebox\myboxA
\newsavebox\myboxB
\newlength\mylenA
\newcommand*\xoverline[2][0.75]{%
    \sbox{\myboxA}{$\m@th#2$}%
    \setbox\myboxB\null
    \ht\myboxB=\ht\myboxA%
    \dp\myboxB=\dp\myboxA%
    \wd\myboxB=#1\wd\myboxA
    \sbox\myboxB{$\m@th\overline{\copy\myboxB}$}
    \setlength\mylenA{\the\wd\myboxA}
    \addtolength\mylenA{-\the\wd\myboxB}%
    \ifdim\wd\myboxB<\wd\myboxA%
       \rlap{\hskip 0.5\mylenA\usebox\myboxB}{\usebox\myboxA}%
    \else
        \hskip -0.5\mylenA\rlap{\usebox\myboxA}{\hskip 0.5\mylenA\usebox\myboxB}%
    \fi}
\makeatother

\DeclareMathOperator{\tr}{tr}
\DeclareMathOperator{\I}{I}
\DeclareMathOperator{\diag}{diag}
\DeclareMathOperator{\Cov}{Cov}
\DeclareMathOperator{\su}{Sum}

\newcommand{\DD}{{\cal D}}
\newcommand{\PP}{{\cal P}}

\newcommand{\eps}{{\varepsilon}}

\newcommand{\smsp}{\hspace{0.3mm}}
\newcommand{\e}{\mathbb{E}}

\newcommand{\p}{\mathbb{P}}

\newcommand{\Reals}{\mathbb{R}}
\newcommand{\Natural}{\mathbb{N}}
\newcommand{\la}{\langle}
\newcommand{\ra}{\rangle}
\newcommand\qed{\hfill\hbox{\rlap{$\sqcap$}$\sqcup$}}

\newcommand{\tsigma}{{\tilde{\sigma}}}

\newcommand{\bigO}[1]{\ensuremath{\mathop{}\mathopen{}\mathcal{O}\mathopen{}\left(#1\right)}}

\newtheorem{lemma}{Lemma}
\newtheorem{theorem}{Theorem}

\theoremstyle{nonumberplain}
\newcommand\specialref{}

\begin{document}

\title{Free energy in the mixed $p$-spin models\\ with vector spins
}
\author{Dmitry Panchenko\thanks{\textsc{\tiny Department of Mathematics, University of Toronto, panchenk@math.toronto.edu. Partially supported by NSERC.}
}\\
}
\date{}
\maketitle
\begin{abstract}
Using the synchronization mechanism developed in the previous work on the Potts spin glass model, we obtain the analogue of the Parisi formula for the free energy in the mixed even $p$-spin models with vector spins, which include the Sherrington-Kirkpatrick model with vector spins interacting through their scalar product. As a special case, this also establishes the sharpness of Talagrand's upper bound for the free energy of multiple mixed $p$-spin systems coupled by constraining their overlaps.
\end{abstract} 
\vspace{0.5cm}

\section{Introduction}

In the previous paper \cite{PPSG}, we computed the free energy in the Potts spin glass model. In this paper, we will extend this result to a more general class of models with vector spins that have arbitrary prior distribution with compact support on $\Reals^\kappa,$ for any $\kappa\geq 1.$ The components of the Hamiltonian on each of the $\kappa$ coordinates of the spin configuration will be given by mixtures of $p$-spin interactions with possibly different sets of inverse temperature parameters. The key step in the computation of the free energy will be exactly the same as in the Potts spin glass, namely, the blocks of overlaps will be forced to synchronize in the infinite-volume limit as a consequence of some special perturbation of the Hamiltonian. This part of the proof will require only cosmetic changes, and we will refer to \cite{PPSG} for the details. Compared to the Potts spin glass, additional difficulties in the general setting are purely technical, mainly due to the fact that we are dealing with arbitrary prior distribution of spins and one has to find the right way to combine techniques from spin glasses and classical large deviations, which takes a little bit of care. 

Let us now describe the model. Fix integer $\kappa\geq 1$ and let $\mu$ be a probability measure on $\Reals^\kappa$ with compact support $\Omega\subseteq \Reals^\kappa.$ A configuration of $N\geq 1$ vector spins will be denoted
\begin{equation}
\sigma = (\sigma_1,\ldots,\sigma_N) \in (\Reals^\kappa)^N,
\end{equation}
the coordinates of each spin $\sigma_i$ will be written as
\begin{equation}
\sigma_i = \bigl(\sigma_i(1),\ldots,\sigma_i(\kappa)\bigr)\in \Reals^\kappa,
\end{equation}
and, for a given $k\leq \kappa$, the configuration of the $k^{\mathrm{th}}$ coordinates will be denoted by
\begin{equation}
\sigma(k)=\bigl(\sigma_1(k),\ldots,\sigma_N(k)\bigr) \in \Reals^N.
\end{equation}
For each $p\geq 2$, let us consider the classical $p$-spin Hamiltonian on $\Reals^N$,
\begin{equation}
H_{N,p}\bigl(\sigma(k)\bigr) = \frac{1}{N^{(p-1)/2}}\sum_{1\leq i_1,\ldots,i_p \leq N} g_{i_1,\ldots,i_p} \sigma_{i_1}(k)\cdots \sigma_{i_p}(k),
\label{VECpurepspin}
\end{equation}
where $(g_{i_1,\ldots,i_p})$ are i.i.d. standard Gaussian for all $p\geq 2$ and $(i_1,\ldots,i_p)$. Notice that these random variables are the same for all coordinates $k\leq \kappa.$ Given sequences $(\beta_p(k))_{p\geq 2}$ for $k\leq \kappa$ of non-negative inverse temperature parameters, we consider mixed $p$-spin Hamiltonians
\begin{equation}
H_{N}^k\bigl(\sigma(k) \bigr) = \sum_{p\geq 2} \beta_p(k) H_{N,p}\bigl(\sigma(k)\bigr).
\label{VECHamHNk}
\end{equation}
We assume that these sequences decrease fast enough to ensure that the series are well defined. For example, if $\Omega\subseteq [-c,c]^\kappa$ then one can take $\beta_p(k)\leq (2c)^{-p}.$ Finally, we define the Hamiltonian of the mixed $p$-spin model with vector spins by
\begin{equation}
H_N(\sigma) = \sum_{k\leq \kappa} H_{N}^k\bigl(\sigma(k) \bigr).
\label{VECVectorHN}
\end{equation}
We will consider only mixed even $p$-spin models, so we will assume that $\beta_p(k)=0$ for all odd $p\geq 3.$ Our main goal will be to compute the limit of the free energy
\begin{equation}
F_N = \frac{1}{N}\smsp \e \log \int_{\Omega^N} \exp H_N(\sigma)\, d\mu^{\otimes N}(\sigma).
\label{VECFEN}
\end{equation}
One can also add a general external field term to the model but, for simplicity of notation, we will omit it. 

\medskip
\noindent
\textbf{Example.}
If one takes $\beta_2(k)=\beta$ for $k\leq \kappa$ and $\beta_p(k)=0$ for $p\geq 3$ then
\begin{equation}
H_N(\sigma) = 
\frac{\beta}{\sqrt{N}}\sum_{1\leq i, j \leq N} g_{ij} (\sigma_i,\sigma_j),
\end{equation}
where $(\sigma_i,\sigma_j)$ is the scalar product of $\sigma_i,\sigma_j\in \Reals^\kappa.$ This choice corresponds to the analogue of the classical Sherrington-Kirkpatrick model \cite{SK} with vector spins interacting through their scalar product. In addition, when the measure $\mu$ is uniform on the standard basis of $\Reals^\kappa$, the model reduces to the Potts spin glass with $\kappa$ orientations, which was considered in \cite{PPSG}. The case when $\mu$ is uniform on the unit circle in $\Reals^2$ is the spin glass analogue of the classical XY or rotor model on the lattice and, when $\mu$ is uniform on the unit sphere in $\Reals^3$, it is the analogue of the classical Heisenberg model on the lattice. The case $\kappa=1$ with the general $\mu$ is the Ghatak-Sherrington model \cite{GS} studied previously in \cite{PGS}.

\medskip
As usual, we will use the upper index, $\sigma^\ell$ for $\ell\geq 1,$ to index sequences of spin configurations. If, for $k,k'\leq \kappa$, we introduce the function
\begin{equation}
\xi_{k,k'}(x) = \sum_{p\geq 2} \beta_p(k) \beta_p(k') x^p
\label{VECxill}
\end{equation}
then it is easy to check that, for two spin configurations $\sigma^\ell$ and $\sigma^{\ell'}$ and for any $k,k'\leq \kappa$,
\begin{equation}
\e H_{N}^k\bigl(\sigma^\ell(k)\bigr) H_{N}^{k'}\bigl(\sigma^{\ell'}(k')\bigr) = N \xi_{k,k'}\bigl(R_{\ell,\ell'}^{k,k'}\bigr),
\label{VECcovVectorHN}
\end{equation}
i.e. the covariance is a function of the overlap between the corresponding coordinates
\begin{equation}
R_{\ell,\ell'}^{k,k'} = \frac{1}{N}\sum_{i\leq N} \sigma_i^{\ell}(k) \sigma_i^{\ell'}(k').
\label{VECRllkk}
\end{equation} 
We will denote the matrix of all such overlaps by $R_{\ell,\ell'}$ or $R(\sigma^\ell,\sigma^{\ell'})$,
\begin{equation}
R_{\ell,\ell'}= R(\sigma^\ell,\sigma^{\ell'}) = \bigl(R_{\ell,\ell'}^{k,k'} \bigr)_{k,k'\leq \kappa}
= \frac{1}{N}\sum_{i\leq N} \sigma_i^{\ell} {\sigma_i^{\ell'}}^T.
\label{VECRll}
\end{equation}
When using matrix operations (transpose, product, etc.) we will always think of vectors as column vectors. In order to state our main result, we need to introduce some notation and definitions.

As in the Potts spin glass in \cite{PPSG}, we will compute the free energy first for a subsystem with constrained self-overlap $R(\sigma,\sigma)$. 
Let us consider the closed convex hull
\begin{equation}
\DD = \mathrm{conv}\Bigl\{  \bigl(\sigma_1(k)\sigma_1(k')\bigr)_{k,k'\leq\kappa} \mid \sigma_1 \in \Omega = \mathrm{supp}(\mu)\Bigr\}
\end{equation}
of $\kappa\times\kappa$ matrices generated by $\sigma_1^{}\sigma_1^T$ for vector spins $\sigma_1\in\Omega$. Clearly, for any $N\geq 1,$ the self-overlap matrix $R(\sigma,\sigma) \in \DD.$ The set $\DD$ is a compact subset of $\Reals^{\kappa\times\kappa}$, as well as the subspace consisting of Gram matrices
\begin{equation}
\Gamma_\kappa = \Bigl\{\gamma \mid \gamma \mbox{ is a $\kappa\times\kappa$ symmetric positive-semidefinite matrix}\Bigr\}. 
\label{VECGammaK}
\end{equation} 
Let $\Pi$ be the space of left-continuous monotone functions on $[0,1]$ with values in $\Gamma_\kappa$,
\begin{equation}
\Pi = \Bigl\{ \pi \colon [0,1]\to \Gamma_\kappa \mid  \pi \mbox{ is left-continuous, } \pi(x)\leq \pi(x') \mbox{ for } x\leq x' \Bigr\},
\end{equation}
where $\pi(x)\leq \pi(x')$ means that $\pi(x')-\pi(x)\in \Gamma_\kappa.$ For $D\in\DD$, we consider
\begin{equation}
\Pi_D = \Bigl\{ \pi\in \Pi \mid \pi(0)=0 \mbox{ and }\, \pi(1) = D \Bigr\}.
\label{VECPid}
\end{equation}
As in \cite{PPSG}, the elements of $\Pi_D$ will play a role of the principle order parameter in the variational formula for the free energy below. We would like to point out that such order parameter already appeared in the physics literature in \cite{FPV}, where a special case of three copies of $p$-spin model was studied in the framework of the Parisi replica method.

A discrete path $\pi\in \Pi_D$ can be encoded by two sequences,
\begin{equation}
x_{-1}=0\leq x_0 \leq \ldots \leq x_{r-1} \leq x_r = 1
\label{VECxs}
\end{equation}
and a monotone sequence of Gram matrices in $\Gamma_\kappa$,
\begin{equation}
0=\gamma_0\leq \gamma_1 \leq \ldots \leq \gamma_{r-1}\leq \gamma_r = D.
\label{VECgammas}
\end{equation}
We can associate to these sequences the path defined by 
\begin{equation}
\mbox{
$\pi(x)= \gamma_j$ for $x_{j-1}<x\leq x_j$ for $0\leq j\leq r,$
}
\label{VECgammapath}
\end{equation}
with $\pi(0)=0$. Recall the function $\xi_{k,k'}$ in (\ref{VECxill}) and denote
\begin{equation}
\theta_{k,k'}(x) = x \xi_{k,k'}'(x)-\xi_{k,k'}(x)= \sum_{p\geq 2} \beta_p(k)\beta_p(k') (p-1)x^p.
\label{VECthetall}
\end{equation}
Given an arbitrary $\kappa\times\kappa$ matrix $A$, we will denote 
\begin{equation}
\xi(A): = \bigl(\xi_{k,k'}(A_{k,k'})\bigr)_{k,k'\leq \kappa},
\label{VECxithetaA}
\end{equation}
and define $\xi'(A)$ and $\theta(A)$ similarly.  If we denote by $\beta_p = (\beta_p(k))_{k\leq \kappa}$ then, for $\gamma\in\Gamma_\kappa$,
\begin{align}
\xi'(\gamma) &= \sum_{p\geq 2}p\gamma^{\,\circ (p-1)}\circ(\beta_p^{} \beta_p^T),
\nonumber
\\
\theta(\gamma) &= \sum_{p\geq 2}(p-1)\gamma^{\,\circ p}\circ(\beta_p^{} \beta_p^T),
\label{VECthetagamma}
\end{align}
where $\circ$ is the Hadamard (element-wise) product and $\gamma^{\,\circ p}$ is the element-wise $p$th power of $\gamma$. An important observation is that these representations imply that the sequences $\xi'(\gamma_j)$ and $\theta(\gamma_j)$ are also non-decreasing in $\Gamma_\kappa$ for $0\leq j\leq r$.

Given a discrete path (\ref{VECgammapath}), let us now consider a sequence of independent Gaussian vectors $z_j = (z_j(k))_{k\leq \kappa}$ for $0\leq j\leq r$ with the covariances
\begin{equation}
\Cov(z_j) = \xi'(\gamma_j) - \xi'(\gamma_{j-1}).
\label{VECzpseq}
\end{equation}
Given $\lambda = (\lambda_{k,k'})_{1\leq k\leq k'\leq \kappa}\in \Reals^{\kappa(\kappa+1)/2}$, let us define 
\begin{equation}
X_r
=
 \log \int_\Omega \exp\Bigl(
\sum_{k\leq\kappa} \sigma_1(k) \sum_{1\leq j \leq r} z_j(k) +\sum_{k\leq k'} \lambda_{k,k'} \sigma_1(k)\sigma_1(k')\Bigr) \,d\mu(\sigma_1)
\label{VECXr}
\end{equation}
and, recursively over $0\leq p\leq r-1,$ define
\begin{equation}
X_j=\frac{1}{x_j}\log \e_j \exp x_j X_{j+1},
\label{VECXp}
\end{equation}
where $\e_j$ denotes the expectation with respect to $z_{j+1}$ only. If $x_j=0$, we interpret this equation  as $X_j = \e_j X_{j+1}.$ Notice that $X_0$ is non-random, and we will denote it by
\begin{equation}
\Phi(\lambda, D, r, x, \gamma) = X_0,
\label{VECPhi}
\end{equation}
making the dependence on all the parameters explicit (the dependence on $D$ here is through the last constraint in (\ref{VECgammas})). For any matrix $A$, we will denote by
\begin{equation}
\su(A) = \sum_{k,k'} A_{k,k'}
\end{equation}
the sum of all its elements. Finally, we define the functional
\begin{equation}
\PP(\lambda,D,r, x, \gamma)
= 
 \Phi(\lambda,D,r, x, \gamma) -\sum_{k\leq k'} \lambda_{k,k'} D_{k,k'}
-\frac{1}{2}
\sum_{0\leq j\leq r-1} x_j \su\bigl(\theta(\gamma_{j+1}) - \theta(\gamma_{j})\bigr).
\label{VECPpar}
\end{equation}
Let us mention right away that, as in the setting of the classical Sherrington-Kirkpatrick model or the Potts spin glass in \cite{PPSG}, one can observe that the functional (\ref{VECPhi}) depends on $(r,x,\gamma)$ only through the path $\pi$ in (\ref{VECgammapath}), so we can denote it by $\Phi(\lambda, D, \pi)$. It was shown in \cite{PPSG} that functionals of this type are Lipschitz with respect to the metric 
\begin{equation}
\Delta(\pi,\pi') = \int_{0}^1 \bigl\|\pi(x)-\pi'(x)\bigr\|_1 \,dx
\label{VECmetric1}
\end{equation}
where $\|\gamma\|_1^{} = \sum_{k,k'} |\gamma_{k,k'}|$. This is a direct analogues of a well-known result of Guerra in \cite{Guerra} (see also \cite{PM} or Theorem 14.11.2 in \cite{SG2}) in the setting of the Sherrington-Kirkpatrick model (see also Lemma \ref{VECLemLC1} below). It was also shown in \cite{PPSG} that a general $\pi\in\Pi_D$ can be discretized in a way that approximates $\pi$ in this metric. Therefore, $\Phi$ can be extended by continuity to all $\pi\in \Pi_D$. Also, rearranging the terms, we can rewrite
\begin{align}
-\sum_{0\leq j\leq r-1} x_j  \su\bigl(\theta(\gamma_{j+1}) - \theta(\gamma_{j})\bigr)
&=
- \su\bigl(\theta(\gamma_{r})\bigr) +  \sum_{1\leq j\leq r} (x_j-x_{j-1}) \su\bigl(\theta(\gamma_{j})\bigr)
\nonumber
\\
&=
-\su\bigl(\theta(D)\bigr) + \int_{0}^1\! \su\bigl(\theta(\pi(x))\bigr)\, dx
\label{VECrearrange}
\end{align}
and, therefore, (\ref{VECPpar}) can be rewritten as
\begin{equation}
\PP(\lambda,D,\pi)
=
 \Phi(\lambda,D,\pi) 
-\sum_{k\leq k'} \lambda_{k,k'} D_{k,k'} 
- \frac{1}{2}\su\bigl(\theta(D)\bigr)
+\frac{1}{2}\int_{0}^1\! \su\bigl(\theta(\pi(x))\bigr) \, dx.
\label{VECPPldpi}
\end{equation}
The following is our main result.
\begin{theorem}\label{VECThFE}
For any $\kappa\geq 1$, the limit of the free energy is given by
\begin{equation}
\lim_{N\to\infty} F_N
= 
\sup_{D\in\DD} \inf_{\lambda,r, x,\gamma} \PP(\lambda,D,r, x, \gamma)
=\sup_{D\in\DD} \inf_{\lambda,\pi\in\Pi_D} \PP(\lambda,D,\pi).
\label{VECParisi}
\end{equation}
\end{theorem}
The formula (\ref{VECParisi}) is the analogue of the classical Parisi formula \cite{Parisi79, Parisi, MPV} for the free energy in the Sherrington-Kirkpatrick model. The upper bound will be a standard application of Guerra's replica symmetry breaking interpolation, and most work will be devoted to the following lower bound. 

Given a subset of spin configurations $S\subseteq \Omega^N$, similarly to (\ref{VECFEN}), we define the free energy constrained to this set of configurations by
\begin{equation}
F_N(S)= \frac{1}{N} \e \log \int_{S} \exp H_{N}(\sigma)\, d\mu^{\otimes N}(\sigma).
\label{VECFES}
\end{equation}
Given $D\in \DD$ and $\eps>0$, we consider an open $\eps$-neighbourhood of $D,$
\begin{equation}
B_\eps(D) = \Bigl\{ \gamma\in\Gamma_\kappa \mid \|\gamma-D\|_\infty <\eps \Bigr\},
\label{VECBepsD}
\end{equation}
with respect to the sup-norm $\|\gamma-D\|_\infty = \max_{k,k'} |\gamma_{k,k'}-D_{k,k'}|$. Let us recall the definition of the overlap matrix in (\ref{VECRll}) and consider the set of spin configurations
\begin{equation}
\Sigma_\eps(D) = \Bigl\{ \sigma\in\Omega^N \mid R(\sigma,\sigma)\in B_\eps(D) \Bigr\}
\label{VECselfoverlapConst}
\end{equation}
with the self-overlap in the $\eps$-neighbourhood of $D$. The lower bound in Theorem \ref{VECThFE} is a direct consequence of the following.
\begin{theorem}\label{VECThFELOW}
For any $D\in\DD$,
\begin{equation}
\lim_{\eps\downarrow 0} \liminf_{N\to\infty} F_N\bigl(\Sigma_\eps(D)\bigr) \geq \inf_{\lambda,r, x,\gamma} \PP(\lambda,D,r, x, \gamma)=\inf_{\lambda,\pi\in\Pi_D} \PP(\lambda,D,\pi).
\label{VECThEqLower}
\end{equation}
\end{theorem}
The proof of the lower bound also works for models with odd $p$-spin interactions, and only the proof of the upper bound uses the convexity of the functions $\xi_{k,k'}$ in (\ref{VECxill}).

Besides the Sherrington-Kirkpatrick model with vector spins mentioned above, perhaps, the most interesting special case included in Theorem \ref{VECThFELOW} corresponds to the uniform measure $\mu$ on $\{-1,+1\}^\kappa$, i.e. multiple copies of the classical mixed $p$-spin model with Ising spins coupled through their overlaps. One of the fundamental ideas in these models is the replica symmetry breaking interpolation invented by Guerra in \cite{Guerra} to show that the Parisi formula \cite{Parisi79, Parisi} is an upper bound on the free energy in the Sherrington-Kirkpatrick model. When Talagrand proved the Parisi formula in \cite{TPF}, the main idea was to generalize Guerra's bound to two copies of the system coupled by fixing their overlap. Since then, various analogues of the Guerra interpolation found many other applications, see for example \cite{MS, Chen-sphere, ChenChaos0, ChenChaos, Chen15, SCLT, FL, PT, PGS, Pspins, PMS, PTChaos, PPSG, T-sphere}. After his seminal work on the Parisi formula, Talagrand proposed in \cite{TalUltra} a generalization of Guerra's bound to multiple copies of the system, possibly at different temperatures, coupled through their overlaps (see Section 15.7 in \cite{SG2}), and suggested a natural approach to other famous problems in spin glasses, such as ultrametricity and chaos, based on looking at the discrepancy between constrained and unconstrained free energies of the system, with constraints violating conjectured properties of the overlaps. However, except for some special cases, this `calculus problem' remains impenetrable. In \cite{TalUltra}, Talagrand raised a possibility that these bounds are not the correct ones, but this possibility is now eliminated by Theorem \ref{VECThFELOW} showing that they are asymptotically sharp. This leaves other possibilities that the analytical structure behind these bounds is extremely non-trivial and remains to be discovered, or that ultrametricity and chaos cannot be detected at the level of the free energy and the probability of spin configurations violating these properties is not exponentially small, although it has been argued in \cite{FPV} (near the critical temperature) and \cite{RC} that both ultrametricity and chaos in temperature can be observed in this way. 

It is interesting to note that, in some sense, we are approaching the sharpness of Talagrand's bounds from the opposite direction, namely, utilizing ultrametricity for the overlaps to study these bounds. Our approach continues the line of ideas originating in another paper of Guerra \cite{GuerraGG}, where the first of the so-called stability properties of the Gibbs measure appeared. The identities for the distribution of the overlaps discovered by Guerra in \cite{GuerraGG} were generalized in \cite{GG} to what are now called the Ghirlanda-Guerra identities. These identities were originally proved on average over temperature, but were later recast by Talagrand in \cite{SG} as a consequence of a small perturbation of the Hamiltonian. This formulation is very powerful because it requires minimal assumptions from the model itself and, as a result, the Ghirlanda-Guerra identities can be proved perturbatively in many other models (the only known example where the Ghirlanda-Guerra identities can be proved non-perturbatively is for generic mixed $p$-spin models, \cite{PGG}). Another related stability property of the Gibbs measure known as the Aizenman-Contucci stochastic stability was discovered in \cite{AC}. The two stability properties can be combined into a unified stability property in the form of the Bolthausen-Sznitman invariance \cite{Bolthausen} in the context of the Ruelle probability cascades, and proved in the context of spin glass models in \cite{ACGG}. The idea of stability turned out to be very fruitful and led to many applications. The first real progress on the ultrametricity problem was made by Arguin and Aizenman in \cite{AA} using the Aizenman-Contucci stochastic stability, under a technical assumption that the overlaps take finitely many values in the infinite-volume limit. A similar result based on the Ghirlanda-Guerra identities was proved in \cite{PGG}, with completely elementary proof discovered later in \cite{PGG2}. The general case turned out to be much harder but it was finally shown in \cite{PUltra} that the Ghirlanda-Guerra identities imply ultrametricity of the overlaps, which means that the Parisi ultrametric ansatz holds perturbatively under minimal assumptions on the model. This led to significant further progress. For example, the Parisi formula was proved in \cite{PPF} for general mixed $p$-spin models including odd $p$-spin interactions, and similar result for the spherical models was proved in \cite{Chen-sphere}. The Ghirlanda-Guerra identities also led to a proof of important symmetries in the setting of diluted spin glass models, namely, the hierarchical exchangeability of pure states, \cite{AP, HEPS}. Combined with a new idea of synchronization of the overlaps, the Parisi ansatz allowed to solve in \cite{PMS} a multi-species version of the Sherrington-Kirkpatrick model introduced in \cite{MS}. The Ghirlanda-Guerra identities played important role in the proof of modified versions of chaos in temperature in \cite{ChenChaos2, ChenChaos}, and a different representation of the Ghirlanda-Guerra identities played a key role in the proof in  \cite{PTChaos} of the first canonical chaos in temperature result for generic mixed even $p$-spin models. A certain overlap-matrix form of the Ghirlanda-Guerra identities was used  to solve the Potts version of the Sherrington-Kirkpatrick model in \cite{PPSG}, and in this paper we will use the same idea to solve the general mixed even $p$-spin models with vector spins, including the sharpness of Talagrand's bound for multiple systems. Hopefully, the observation that Talagrand's bounds are sharp will serve as a further motivation to try to understand if they can be exploited to study mixed $p$-spin models non-perturbatively and, for example, prove chaos in temperature for all mixed $p$-spin models.

As we mentioned above, the main idea of the proof is identical to the setting of the Potts spin glass \cite{PPSG}, and the corresponding parts of the proof will be only recalled briefly or sketched, especially, when they are slightly modified. The main new technical difficulty comes from the fact that, for a general measure $\mu$, we can constrain the self-overlap as in (\ref{VECselfoverlapConst}) only up to some small $\eps>0$, while the covariance structure of various cavity fields in the usual cavity computations must be constrained more precisely in the limit, in order for spin glass techniques to work. Once we start cavity computations in Section \ref{VECSec3label}, we will explain these issues in more detail to motivate the sections that follow. In fact, we will break the cavity computations of the lower bound in three sections, Sections \ref{VECSec3label},  \ref{VECSec6label} and  \ref{VECSec8label}, which will alternate with necessary technical results proved in between. In Section \ref{VECSec4label}, we will construct a certain modification of the spin configurations designed to make the main idea work smoothly in the present setting and, in Section \ref{VECSec5label}, we will reformulate the perturbation and synchronization mechanisms developed in the setting of the Potts spin glass. Section \ref{VECSec7label} will be devoted to some standard large deviation techniques for the functionals that appear in the infinite-volume limit. We begin in Section \ref{VECSec2label} with the analogue of Guerra's replica symmetry breaking interpolation and the proof of the upper bound. 

\medskip
\noindent
\textbf{Acknowledgements.} The author would like to thank Giorgio Parisi for several comments which led to improvement of the paper.

\section{Upper bound via Guerra's interpolation}\label{VECSec2label}

\noindent
\textbf{Remark.} Throughout the paper, we will denote by $L$ any constant that does not depend on any individual parameters, such as $D\in \DD$ or $N$, but depends only on the global parameters of the model, such as the dimension $\kappa,$ the covariance structure of the Hamiltonian and the size of the support of the measure $\mu$. The constant can change even within the same equation.

\medskip
The proof of the upper bound is, essentially, identical to Section 15.7 in \cite{SG2}. By continuity, in the rest of the paper we will assume that the inequalities in (\ref{VECxs}) are strict,
\begin{equation}
x_{-1}=0< x_0 < \ldots < x_{r-1} < x_r = 1.
\label{VECxsstrict}
\end{equation}
Let $(v_\alpha)_{\alpha\in \Natural^r}$ be the weights of the Ruelle probability cascades \cite{Ruelle} corresponding to the sequence (\ref{VECxsstrict}) (see e.g. Section 2.3 in \cite{SKmodel} for the definition). For $\alpha^1, \alpha^2\in \Natural^r$, we denote
\begin{equation}
\alpha^1\wedge \alpha^2 = \min\Bigl\{0\leq j \leq r  \mid  \alpha_1^1= \alpha_1^2, \ldots, \alpha_{j}^1 = \alpha_{j}^2, \alpha_{j+1}^1 \not = \alpha_{j+1}^2 \Bigr\},
\end{equation}
where $\alpha^1\wedge \alpha^2 =r$ if $\alpha^1=\alpha^2$. We observed in (\ref{VECthetagamma}) that the sequences $\xi'(\gamma_j)$ and $\theta(\gamma_j)$ are non-decreasing in $\Gamma_\kappa$ for $0\leq j\leq r$. As a result, there exist Gaussian processes 
\begin{equation}
Z(\alpha)=\bigl(Z^k(\alpha)\bigr)_{k\leq\kappa} \,\mbox{ and }\, Y(\alpha),
\end{equation}
both indexed by $\alpha\in\Natural^r$, with the covariances
\begin{align}
\Cov\bigl(Z(\alpha^1), Z(\alpha^2)\bigr) &= \xi'\bigl(\gamma_{\alpha^1\wedge\alpha^2}\bigr),
\nonumber
\\
\Cov\bigl(Y(\alpha^1), Y(\alpha^2)\bigr) &= \su\bigl(\theta(\gamma_{\alpha^1\wedge\alpha^2}) \bigr).
\label{VECCD}
\end{align}
Let $Z_i(\alpha)$ be independent copies of the process $Z(\alpha)$, also independent of $Y(\alpha)$. For $0\leq t\leq 1$, consider an interpolating Hamiltonian defined on $\Omega^N\times \Natural^r$ by
\begin{equation}
H_{N,t}(\sigma,\alpha) = 
\sqrt{t} H_N(\sigma) + \sqrt{1-t} \sum_{i\leq N} 
\sum_{k\leq \kappa} \sigma_i(k) Z_{i}^k(\alpha) +\sqrt{t} \sqrt{N} Y(\alpha).
\label{VECHNt}
\end{equation}
Similarly to (\ref{VECFES}), we define the interpolating free energy constrained to the set $S\subseteq \Omega^N$, 
\begin{equation}
\varphi_S(t)= \frac{1}{N} \e \log \sum_{\alpha\in\Natural^r} v_\alpha \int_{S} \exp H_{N,t}(\sigma,\alpha)\, d\mu^{\otimes N}(\sigma).
\label{VECFEDt}
\end{equation}
Recall the definition of the set $\Sigma_\eps(D)$ in (\ref{VECselfoverlapConst}). We begin with the following.
\begin{lemma}\label{VECLemGUP}
The derivative of the function $\varphi(t)$ in (\ref{VECFEDt}) with $S=\Sigma_\eps(D)$ satisfies $\varphi'(t)\leq L\eps$.
\end{lemma}
\textbf{Proof.} Let us denote by $\la\, \cdot\,\ra_t$ the average with respect to the measure 
$$
G_t(d\sigma,\alpha) \sim v_\alpha \exp H_{N, t}(\sigma,\alpha)\,d\mu^{\otimes N}(\sigma).
$$
on $\Sigma_\eps(D)\times \Natural^r$. Then, for $0<t<1$,
$$
\varphi'(t) = \frac{1}{N}\e \Bigl\la \frac{\partial H_{N, t}(\sigma,\alpha)}{\partial t} \Bigr\ra_t.
$$
From the definition of $H_N(\sigma)$ in (\ref{VECVectorHN}) and (\ref{VECcovVectorHN}),
\begin{equation}
\e H_N(\sigma^1) H_N(\sigma^2) = N\su\bigl(\xi(R_{1,2})\bigr).
\label{VECCovHNsum}
\end{equation}
Similarly, from the definition (\ref{VECCD}),
\begin{equation}
\e \sum_{i\leq N} \sum_{k\leq \kappa} \sigma_i^1(k) Z_{i}^k(\alpha^1)
\sum_{i\leq N} \sum_{k\leq \kappa} \sigma_i^2(k) Z_{i}^k(\alpha^2)
=
N\su \bigl( R_{1,2}\circ \xi'(\gamma_{\alpha^1\wedge\alpha^2})\bigr).
\end{equation}
Using these equations and recalling the covariance of $Y(\alpha)$ in (\ref{VECCD}),
$$
\frac{1}{N}\, \e \frac{\partial H_{N, t}(\sigma^1,\alpha^1)}{\partial t} H_{N, t}(\sigma^2,\alpha^2)
=
\frac{1}{2} \su\Bigl(\xi(R_{1,2}) - R_{1,2}\circ \xi'(\gamma_{\alpha^1\wedge\alpha^2}) + \theta(\gamma_{\alpha^1\wedge\alpha^2}) \Bigr).
$$
By the usual Gaussian integration by parts (see e.g. Lemma 1.1 in \cite{SKmodel}),
\begin{align*}
\varphi'(t) = &\,\,
\frac{1}{2}
\e \Bigl\la \su\Bigl(\xi(R_{1,1}) - R_{1,1}\circ \xi'(\gamma_{\alpha^1\wedge\alpha^1}) + \theta(\gamma_{\alpha^1\wedge\alpha^1}) \Bigr)
\Bigr\ra_t
\\
& - \frac{1}{2}
\e \Bigl\la \su\Bigl(\xi(R_{1,2}) - R_{1,2}\circ \xi'(\gamma_{\alpha^1\wedge\alpha^2}) + \theta(\gamma_{\alpha^1\wedge\alpha^2}) \Bigr)
\Bigr\ra_t.
\end{align*}
Since $\theta_{k,k'}(x) =x\xi_{k,k'}'(x)-\xi_{k,k'}(x)$ for all $k,k'\leq \kappa$, $\gamma_{\alpha^1\wedge\alpha^1} = \gamma_r = D$ and $R_{1,1} = R(\sigma^1,\sigma^1)\in B_\eps(D)$ for $\sigma^1\in \Sigma_\eps(D)$, the first term is bounded by $L\eps$. We also have $\xi_{k,k'}(a)-a\xi_{k,k'}'(b)+\theta_{k,k'}(b)\geq 0$ by convexity of $\xi_{k,k'}$, so the second term is negative and this finishes the proof.
\qed

\begin{lemma}\label{VECLemAlmostUp}
For any $\lambda = (\lambda_{k,k'})_{1\leq k\leq k'\leq \kappa}\in \Reals^{\kappa(\kappa+1)/2}$,
\begin{equation}
F_N\bigl(\Sigma_\eps(D)\bigr) \leq  L\eps + \eps \|\lambda\|_1 + \PP(\lambda,D,r, x, \gamma).
\label{VECalmostUP}
\end{equation}
\end{lemma}
\textbf{Proof.}
At the beginning of the interpolation at $t=1$,
\begin{equation}
\varphi_{\Sigma_\eps(D)}(1)
=
F_N\bigl(\Sigma_\eps(D)\bigr)
+
\frac{1}{N}\smsp \e\log \sum_{\alpha\in\Natural^r} v_{\alpha} 
\exp \sqrt{N} Y(\alpha).
\label{VECfunp20}
\end{equation}
The standard properties of the Ruelle probability cascades (see Section 2.3 and the proof of Lemma 3.1 in \cite{SKmodel}) together with the covariance structure (\ref{VECCD}) imply that 
\begin{equation}
\frac{1}{N}\smsp \e\log \sum_{\alpha\in\Natural^r} v_{\alpha} 
\exp \sqrt{N} Y(\alpha)
=
\frac{1}{2}
\sum_{0\leq j\leq r-1} x_j  \su\bigl( \theta(\gamma_{j+1}) - \theta(\gamma_{j})\bigr).
\label{VECfunp2}
\end{equation}
Next, let us consider
$$
\varphi_{\Sigma_\eps(D)}(0)
=
\frac{1}{N} \e \log \sum_{\alpha\in\Natural^r} v_\alpha \int_{\Sigma_\eps(D)} \exp \Bigl( \sum_{i\leq N} 
\sum_{k\leq \kappa} \sigma_i(k) Z_{i}^k(\alpha) \Bigr) \,d\mu^{\otimes N}(\sigma).
$$
For any $\sigma\in \Sigma_\eps(D)$ and $\lambda \in \Reals^{\kappa(\kappa+1)/2}$,
$$
-\sum_{k\leq k'} \lambda_{k,k'} D_{k,k'}
+ \frac{1}{N}\sum_{i\leq N} \sum_{k\leq k'} \lambda_{k,k'} \sigma_i(k)\sigma_i(k')
\geq -\eps \|\lambda\|_1
$$
and, therefore,
\begin{align*}
\varphi_{\Sigma_\eps(D)}(0)
\leq &\,\,\,
\eps \|\lambda\|_1 -\sum_{k\leq k'} \lambda_{k,k'} D_{k,k'} 
\\
&
+
\frac{1}{N} \e \log \sum_{\alpha\in\Natural^r} v_\alpha \int_{\Omega^N} \exp \sum_{i\leq N}
\Bigl( \sum_{k\leq \kappa} \sigma_i(k) Z_{i}^k(\alpha) + \sum_{k\leq k'} \lambda_{k,k'} \sigma_i(k)\sigma_i(k') \Bigr)\, d\mu^{\otimes N}(\sigma).
\end{align*}
If we introduce the notation
$$
X_i^\alpha =
\int_{\Omega} \exp
\Bigl( \sum_{k\leq \kappa} \sigma_i(k) Z_{i}^k(\alpha) + \sum_{k\leq k'} \lambda_{k,k'} \sigma_i(k)\sigma_i(k') \Bigr)\, d\mu(\sigma_i)
$$
then this upper bound can be rewritten as
$$
\varphi_{\Sigma_\eps(D)}(0)
\leq 
\eps \|\lambda\|_1 -\sum_{k\leq k'} \lambda_{k,k'} D_{k,k'}  +
\frac{1}{N} \e \log \sum_{\alpha\in\Natural^r} v_\alpha 
\prod_{i\leq N} X_{i}^\alpha.
$$
Standard properties of the Ruelle probability cascades (see Section 2.3 in \cite{SKmodel}) imply that
$$
\frac{1}{N} \e \log \sum_{\alpha\in\Natural^r} v_\alpha \prod_{i\leq N} X_{i}^\alpha
=
\e \log \sum_{\alpha\in\Natural^r} v_\alpha X_{1}^\alpha = X_0,
$$
where $X_0 = \Phi(\lambda, D, r, x, \gamma)$ was defined in (\ref{VECPhi}) and, therefore,
\begin{equation}
\varphi_{\Sigma_\eps(D)}(0)
\leq 
\eps \|\lambda\|_1 -\sum_{k\leq k'} \lambda_{k,k'} D_{k,k'} 
+ \Phi(\lambda, D, r, x, \gamma).
\label{VECphi0bound}
\end{equation}
Together with (\ref{VECfunp20}), (\ref{VECfunp2}) and Lemma \ref{VECLemGUP} this implies that $F_N(\Sigma_\eps(D))$ is bounded by
$$
L\eps + \eps \|\lambda\|_1 -\sum_{k\leq k'} \lambda_{k,k'} D_{k,k'} 
+ \Phi(\lambda, D, r, x, \gamma) 
- \frac{1}{2} \sum_{0\leq j\leq r-1} x_j  \su\bigl( \theta(\gamma_{j+1}) - \theta(\gamma_{j})\bigr),
$$
which finishes the proof.
\qed

\medskip
We are now ready to prove the upper bound in Theorem \ref{VECThFE}.
\begin{lemma}\label{VECThFEUP}
For any $\kappa\geq 1$, the free energy satisfies
\begin{equation}
\limsup_{N\to\infty} F_N
\leq
\sup_{D} \inf_{\lambda,r, x,\gamma} \PP(\lambda,D,r, x, \gamma).
\label{VECParisiUP}
\end{equation}
\end{lemma}
\textbf{Proof.}
Fix $\delta>0$ and, for $D\in \DD$, let
$$
\PP_\delta(D) = \max\Bigl(-\frac{1}{\delta}, \delta + \inf_{\lambda,r, x,\gamma} \PP(\lambda,D,r, x, \gamma)\Bigr).
$$
For each $D\in \DD$, one can find $\lambda_D, r_D, x_D, \gamma_D$ such that $\PP(\lambda_D, D, r_D, x_D, \gamma_D)\leq \PP_\delta(D).$ If $L$ is a constant in (\ref{VECalmostUP}), let $\eps_D>0$ be such that $\eps_D(L+\|\lambda_D\|_1)\leq \delta.$ Lemma \ref{VECLemAlmostUp} then implies that 
$$
F_N\bigl(\Sigma_{\eps_D}(D)\bigr)\leq \delta + \PP_\delta(D).
$$
Since the collection of sets $B_{\eps_D}(D)$ for $D\in\DD$ form an open cover of $\DD$ and $\DD$ is compact, we can find a finite subcover indexed by some $D_1,\ldots, D_n \in \DD.$ Consider the random  free energy with spin configurations constrained to the set $S$,
$$
\tilde{F}_N(S)= \frac{1}{N} \log \int_{S} \exp H_{N}(\sigma)\, d\mu^{\otimes N}(\sigma).
$$
Since the union of $\Sigma_{\eps_{D_i}}(D_i)$ for $i\leq n$ covers $\Omega^N$,
$$
\tilde{F}_N(\Omega^N)
\leq 
\frac{\log n}{N} + \max_{i\leq n} \tilde{F}_N\bigl(\Sigma_{\eps_{D_i}}(D_i)\bigr).
$$
By the Gaussian concentration inequalities, $\tilde{F}_N(S)$ deviates from its expectation $F_N(S)$ by more than $1/\sqrt{N}$ with exponentially small probability of the order $e^{-N/L}$, where the constant $L$ does not depend on the set $S$. With the above inequalities, this implies that
$$
{F}_N
\leq 
\frac{2}{\sqrt{N}} + \frac{\log n}{N} + \max_{i\leq n} {F}_N\bigl(\Sigma_{\eps_{D_i}}(D_i)\bigr)
\leq 
\frac{2}{\sqrt{N}} + \frac{\log n}{N} 
+ \delta + \max_{i\leq n} \PP_\delta(D_i).
$$
Therefore, $$\limsup_{N\to\infty} F_N \leq \delta + \sup_D \PP_\delta(D)$$ and letting $\delta\downarrow 0$ finishes the proof.
\qed

\section{Cavity computation, part 1}\label{VECSec3label}

The proof of the lower bound in Theorem \ref{VECThFELOW} will take up the rest of the paper, and we will start with a standard Aizenman-Sims-Starr cavity computation \cite{AS2} in the form that appeared, for example, in \cite{SKcoupled, Chen-sphere, PPSG}. Let us make the dependence of $\Sigma_\eps(D)$ in (\ref{VECselfoverlapConst}) on $N$ explicit, $\Sigma^N_\eps(D)$, and denote
\begin{equation}
Z_{N}(\eps,D) = \int_{\Sigma^N_\eps(D)} \exp H_N(\sigma)\, d\mu^{\otimes N}(\sigma),
\label{VECZpartE}
\end{equation}
so that $F_N(\Sigma^N_\eps(D)) = N^{-1}\e\log Z_N(\eps, D).$ We start with an obvious inequality,
\begin{equation}
\liminf_{N\to\infty} F_N\bigl(\Sigma^N_\eps(D)\bigr) \geq \frac{1}{M}\liminf_{N\to\infty} \Bigl(\e \log Z_{N+M}(\eps, D) - \e \log Z_{N}(\eps, D)\Bigr),
\label{VECFNAN}
\end{equation}
where $M$ on the right hand side is fixed. Let us write spin configurations in $\Omega^{N+M}$ as $\rho=(\sigma,\tau)$ for $\sigma\in\Omega^N$ and $\tau\in\Omega^M.$ Using that
$$
R(\rho,\rho) = \frac{N}{N+M}R(\sigma,\sigma)+\frac{M}{N+M}R(\tau,\tau),
$$
we get that
$$
\Bigl\{\rho \mid R(\rho,\rho)\in B_\eps(D) \Bigr\}
\supseteq 
\Bigl\{\sigma \mid R(\sigma,\sigma)\in B_\eps(D) \Bigr\}
\times
\Bigl\{\tau \mid R(\tau,\tau)\in B_\eps(D) \Bigr\}
$$
and, therefore,
$$
Z_{N+M}(\eps,D) \geq \int_{\Sigma^N_\eps(D)} \int_{\Sigma^M_\eps(D)} 
\exp H_{N+M}(\sigma,\tau)\, d\mu^{\otimes M}(\tau) d\mu^{\otimes N}(\sigma).
$$
This allows to decrease the lower bound in (\ref{VECFNAN}) to
\begin{align}
\liminf_{N\to\infty} F_N(\Sigma^N_\eps(D))
\geq
\liminf_{N\to\infty} \frac{1}{M} \Bigl(
&
\e \log \int_{\Sigma^N_\eps(D)} \int_{\Sigma^M_\eps(D)} 
\exp H_{N+M}(\sigma,\tau)\, d\mu^{\otimes M}(\tau) d\mu^{\otimes N}(\sigma)
\nonumber
\\
&
- \e \log \int_{\Sigma^N_\eps(D)} \exp H_{N}(\sigma)\, d\mu^{\otimes N}(\sigma)
\Bigr).
\label{VECFMNAN}
\end{align}
Then one can do the usual calculation as in the Aizenman-Sim-Starr representation \cite{AS2} (see e.g. Section 1.3 in \cite{SKmodel}), separating the Hamiltonian
\begin{equation}
H_{N+M}(\sigma,\tau) = H_N'(\sigma) + \sum_{i\leq M}\sum_{k\leq\kappa} \tau_i(k) Z_{i}^k(\sigma) + r(\tau)
\label{VECcommonH}
\end{equation}
into three types of terms -- that depend only on $\sigma$, the ones where only one spin $\tau_i$ appears, and the ones where more than two coordinates of $\tau$ appear. Of course, $Z_{i}^k(\sigma)$ here depends only on the $k$th coordinate $\sigma(k)$ of the configuration $\sigma$, but the dependence on $k$ is already reflected in the upper index (this includes the dependence on the parameters $(\beta_p(k))_{p\geq 2}$ in (\ref{VECHamHNk})). The term $r(\tau)$ can be omitted because it is of a small order as $N\to\infty$. The Gaussian process $H_N'(\sigma)$ is defined just like $H_N(\sigma)$, only with scalings in (\ref{VECpurepspin}) by the powers of $N+M$ instead of $N$. As a result, one can decompose (in distribution),
\begin{equation}
H_N(\sigma) \stackrel{d}{=}
H_N'(\sigma) + \sqrt{M} \sum_{k\leq \kappa} Y^k(\sigma),
\label{VECcommonH2}
\end{equation}
for some Gaussian processes $Y^k(\sigma)$ independent of $H_N'(\sigma).$ One can easily check (see e.g. Section 3.5 in \cite{SKmodel} for a similar computation) that, for $k,k'\leq \kappa$,
\begin{align}
\e Z_i^k(\sigma^{\ell}) Z_i^{k'}(\sigma^{\ell'})
&= 
\xi_{k,k'}'(R_{\ell,\ell'}^{k,k'}) + \bigO{\frac{M}{N}},
\label{VECCovz}
\\
\e Y^k(\sigma^\ell) Y^{k'}(\sigma^{\ell'}) 
&= 
\theta_{k,k'}(R_{\ell,\ell'}^{k,k'})
+ \bigO{\frac{M}{N}}.
\label{VECCovy}
\end{align}
If we define $Y(\sigma) = \sum_{k\leq \kappa} Y^k(\sigma)$ then
\begin{equation}
\e Y(\sigma^\ell) Y(\sigma^{\ell'}) 
= 
\su\bigl(\theta(R_{\ell,\ell'})\bigr)
+ \bigO{\frac{M}{N}}.
\label{VECCovysum}
\end{equation}
One can redefine the processes $Z_i^k$ and $Y^k$ to have the covariances without the error terms $\bigO{M/N}$, since this does not affect the right hand side of (\ref{VECFMNAN}), which we assume from now on.

Consider the Gibbs measure on $\Sigma^N_\eps(D)$ corresponding to the Hamiltonian $H_N'(\sigma)$ in (\ref{VECcommonH}),
\begin{equation}
dG_N(\sigma) = \frac{\exp H_N'(\sigma)d\mu^{\otimes N}(\sigma)}{Z_N'(\eps, D)},
\,\mbox{ where }\,
Z_N'(\eps, D) = \int_{\Sigma^{N}_\eps(D)} \exp H_N'(\sigma)\, d\mu^{\otimes N}(\sigma)
\label{VECMeasureGNprime}
\end{equation}
and let us denote by $\la\,\cdot\,\ra_N$ the average with respect to $G_N$. Using representations (\ref{VECcommonH}) and (\ref{VECcommonH2}) (omitting the negligible term $r(\tau)$) and dividing inside both $\log$arithms by $Z_N'(\eps, D)$, we can rewrite the quantity on the right hand side of (\ref{VECFMNAN}) as
\begin{equation}
\frac{1}{M}
\Bigl(
\e \log \Bigl\la
\int_{\Sigma^{M}_\eps(D)} \exp\Bigl( \sum_{i\leq M}\sum_{k\leq\kappa} \tau_i(k) Z_{i}^k(\sigma)
\Bigr)\, d\mu^{\otimes M}(\tau)
\Bigr\ra_N
-
\e \log \Bigl\la \exp \sqrt{M} Y(\sigma) \Bigr\ra_N
\Bigr).
\label{VECAS2repr}
\end{equation}
Both terms here are continuous functionals of the distribution of the overlap array $(R_{\ell,\ell'})_{\ell,\ell'\geq 1}$ under the measure $\e (G_N)^{\otimes \infty}$ (see e.g. the proof of Theorem 1.3 in \cite{SKmodel}), so in order to understand the limit $N\to\infty$, we need to understand the behaviour of this distribution. This will be achieved via the main idea used to solve the Potts spin glass in \cite{PPSG}, namely, a special perturbation of the Hamiltonian $H_N'(\sigma)$ which will ensure the validity of the overlap-matrix version of the Ghirlanda-Guerra identities. 

However, there is an issue we have to deal with that did not arise in \cite{PPSG}. Namely, the diagonal overlap blocks $R_{\ell,\ell}=R(\sigma^\ell,\sigma^\ell)$ for replicas $\sigma^\ell$ sampled from the measure $G_N$ are now not fixed, since we only constrain them to be in the $\eps$-neighbourhood $B_\eps(D)$ of $D$, and they can not satisfy the Ghirlanda-Guerra identities that are central to the whole argument. We will resolve this issue by mapping configurations $\sigma\in\Sigma_\eps^N(D)$ into configurations $\tsigma$ such that $R(\tsigma,\tsigma)$ is fixed. We need to do this in a way that controls global distortion and does not affect the overlaps $R(\sigma^\ell,\sigma^{\ell'})$ much. Once we see how this can be done, the processes $Z_i^k(\sigma)$ and $Y(\sigma)$ in (\ref{VECAS2repr}) will be replaced by $Z_i^k(\tsigma)$ and $Y(\tsigma)$ with the covariance depending on the overlaps $R(\tsigma^\ell,\tsigma^{\ell'})$. In particular, since the Ghirlanda-Guerra identities is a property of the perturbation of the Hamiltonian, this perturbation will need to be directly defined in terms of $\tsigma.$

This introduces another issue we have to be aware of when we define the mapping $\sigma\to\tsigma.$ As in the Potts spin glass \cite{PPSG}, in the above cavity computation, the Hamiltonian $H_{N+M}(\sigma,\tau)$ will have a perturbation term $s_{N+M} h_{N+M}(\tilde{\rho})$ with $\rho=(\sigma,\tau)$, while the Hamiltonian $H_{N}(\sigma)$ will come with the perturbation term $s_{N} h_{N}(\tsigma)$ and, as usual, in the first term in (\ref{VECFNAN}) we will replace $s_{N+M} h_{N+M}(\tilde{\rho})$ by $s_{N} h_{N}(\tsigma)$. Since we will take $s_N = N^\gamma$ for any $1/4<\gamma<1/2$, which is not small, and the covariance of $h_N$ will be a continuous function of the overlap $R(\tsigma^\ell,\tsigma^{\ell'})$, in order to make this work, we will need the difference between $R(\tsigma^\ell,\tsigma^{\ell'})$ and $R(\tilde{\rho}^\ell,\tilde{\rho}^{\ell'})$ to be of the order $1/N.$ Since the difference between $R(\sigma^\ell,\sigma^{\ell'})$ and $R({\rho}^\ell,{\rho}^{\ell'})$ is of order $1/N$, this again amounts to controlling the distortion of the map $\sigma\to\tsigma.$ We will come back to the cavity computation after we resolve these issues and recall the matrix Ghirlanda-Guerra identities from \cite{PPSG}.

\section{Modification of spin configurations}\label{VECSec4label}

Given a matrix $D\in\DD\subseteq \Gamma_\kappa$, let
\begin{equation}
D=Q\Lambda Q^T,\,\,
\Lambda = \diag(\lambda_1,\ldots,\lambda_\kappa),
\label{VECQLQ}
\end{equation}
be its eigendecomposition. Without loss of generality, suppose that the eigenvalues are arranged in the decreasing order, $\lambda_1\geq \ldots\geq \lambda_\kappa,$ and, given $\eps>0$, let $0\leq m\leq \kappa$ be such that $\lambda_m\geq \sqrt{\eps}$ and $\lambda_{m+1}<\sqrt{\eps}$. Let us define
\begin{equation}
D_\eps=Q\Lambda_\eps Q^T,\,\,
\Lambda = \diag(\lambda_1,\ldots,\lambda_m,0,\ldots,0).
\label{VECQLQeps}
\end{equation}
Given any $\sigma\in\Sigma_\eps(D)$, which means that self-overlap $R(\sigma,\sigma)\in B_\eps(D)$, we will construct a $\kappa\times \kappa$ matrix $A$ such that the self-overlap of $A\sigma =(A\sigma_i)_{i\leq N}$ is equal to $D_\eps$,
\begin{equation}
R(A\sigma, A\sigma) = \frac{1}{N}\sum_{i\leq N} (A\sigma_i)(A\sigma_i)^T =
A R(\sigma,\sigma) A^T= D_\eps,
\end{equation}
and such that $A$ has small distortion in the sense explained below. The reason we removed the eigenvalues smaller than $\sqrt{\eps}$ in $D$ is precisely to ensure that $A$ has small distortion. These small eigenvalues will be reintroduced at the very end of the computation of the lower bound, using continuity properties of the functionals involved. The matrix $A$ will depend on $\sigma$ only through the self-overlap $R(\sigma,\sigma)$, and we will denote $A$ by 
$$
A_\sigma \,\mbox{ or }\, A(R(\sigma,\sigma))
$$ 
when we need to make this dependence explicit.

First of all, small distortion means that the overlaps of $\sigma$ with other configurations in $\Omega^N$ should not change much when $\sigma$ is replaced by $\tsigma =( A_\sigma \sigma_i)_{i\leq N}.$ A convenient way to control the difference is as follows. If $\rho\in\Omega^N$ and $v = \tsigma - \sigma$ then
\begin{align}
\Bigl\| \frac{1}{N}\sum_{i\leq N} \tsigma_i^{} \rho_i^T - \frac{1}{N}\sum_{i\leq N} \sigma_i^{}  \rho_i^T \Bigr\|_{HS}
& =
\Bigl\| \frac{1}{N}\sum_{i\leq N} v_i^{}  \rho_i^T \Bigr\|_{HS}
\leq
\frac{1}{N}\sum_{i\leq N} \bigl\| v_i^{}  \rho_i^T \bigr\|_{HS}
\label{VECoverchange}
\\
&=
\frac{1}{N}\sum_{i\leq N} \| v_i \| \| \rho_i \|
\leq 
\frac{L}{N}\sum_{i\leq N} \| v_i \|
\nonumber
\\
&\leq 
L\Bigl(\frac{1}{N}\sum_{i\leq N} \| v_i \|^2 \Bigr)^{1/2}
=
L\Bigl(\tr\Bigl( \frac{1}{N}\sum_{i\leq N} v_i^{} v_i^T  \Bigr) \Bigr)^{1/2}.
\nonumber
\end{align}
Since $v = \tsigma - \sigma = (A_\sigma-I)\sigma$,
\begin{equation}
\frac{1}{N}\sum_{i\leq N} v_i^{} v_i^T 
=
(A_\sigma-I) R(\sigma,\sigma) (A_\sigma-I)^T,
\end{equation}
and we can control the difference of the overlaps via the trace of this matrix. Another piece of information about the map $A_\sigma$ that we will need is motivated by the following question. Suppose that we have two pairs of configurations $\sigma^1,\sigma^2$ and $\rho^1,\rho^2$ that are close to each other in the sense that their overlaps $R(\sigma^1,\sigma^2)$ and $R(\rho^1,\rho^2)$ and self-overlaps $R(\sigma^j,\sigma^j)$ and $R(\rho^j,\rho^j)$ are close to each other. Then, how close will the overlaps 
\begin{equation}
R(A_{\sigma^1}\sigma^1,A_{\sigma^2}\sigma^2) = A_{\sigma^1}^{} R(\sigma^1,\sigma^2) A_{\sigma^2}^T
\,\mbox{ and }\,
R(A_{\rho^1}\rho^1,A_{\rho^2}\rho^2) = A_{\rho^1}^{} R(\rho^1,\rho^2) A_{\rho^2}^T
\label{VECRoverR}
\end{equation}
be after we apply the corresponding transformations to all the configurations? For this, we will need to control the sup-norms $\|A_{\sigma^1}-A_{\sigma^2}\|_\infty$, which will be bounded in terms of the sup-norm $\|R(\sigma^1,\sigma^1)-R(\sigma^2,\sigma^2)\|_\infty$.
\begin{lemma}\label{VECLemModifyA}
For each $R\in B_\eps(D)$ there exists a matrix $A=A(R)\in\Gamma_\kappa$ such that $A R A^T= D_\eps,$
\begin{equation}
\tr\bigl(
(A-I) \,R\,(A-I)^T
\bigr)\leq L\sqrt{\eps}
\label{VECdiscone}
\end{equation}
and, for any $R_1, R_2\in B_\eps(D)$,
\begin{equation}
\|A(R_1)-A(R_2)\|_\infty
\leq
\frac{L}{\eps}
\|R_1-R_2\|_\infty.
\label{VECdisctwo}
\end{equation}
\end{lemma}
\textbf{Proof.}
Recall the decomposition in (\ref{VECQLQ}). Let us change the coordinate system by applying the transformation $Q^T R Q$ to all matrices, which does not change the trace and changes the sup norm $\|R\|_\infty$ only up to a constant factor. In particular, $Q^T R Q\in B_{L\eps}(\Lambda)$. Once we define  $A(Q^T R Q)$, we can go back and define $A(R)=QA(Q^T R Q)Q^T.$ As a result, from now on we assume that $D=\Lambda$ and $R\in B_\eps(\Lambda).$ 

Let us recall (\ref{VECQLQeps}) and let us denote $\Lambda_m = \diag(\lambda_1,\ldots,\lambda_m).$ If $Q=(R_{k,k'})_{k,k'\leq m}$ is the matrix consisting of the first $m$ rows and columns of $R$, then $Q\in B_\eps(\Lambda_m).$ Suppose we can find $m\times m$ matrix $B=B(Q)\in \Gamma_m$ such that $B Q B^T= \Lambda_m,$
\begin{equation}
\tr\bigl(
(B-I) \,Q\,(B-I)^T
\bigr)\leq L\sqrt{\eps}
\label{VECdisconetwo}
\end{equation}
and, for any $Q_1, Q_2\in B_\eps(\Lambda_m)$, 
\begin{equation}
\|B(Q_1)-B(Q_2)\|_\infty
\leq
\frac{L}{\eps}
\|Q_1-Q_2\|_\infty.
\label{VECdisctwotwo}
\end{equation}
Then, we will define $A(R)$ by extending $B(Q)$ by all zeros in rows and columns from $m+1$ to $\kappa.$ Then (\ref{VECdisctwo}) will, obviously, follow from (\ref{VECdisctwotwo}). As for (\ref{VECdiscone}), if we denote by $T=(R_{k,k'})_{k,k'\geq m+1}$ the matrix consisting of the last $\kappa-m$ rows and columns of $R$, then
$$
\tr\bigl(
(A-I) \,R\,(A-I)^T
\bigr)
=
\tr\bigl(
(B-I) \,Q\,(B-I)^T
\bigr)
+\tr(T).
$$
However, since $R\in B_\eps(\Lambda),$ we have
$$
\tr(T)\leq (\kappa-m)\eps + \lambda_{m+1}+\ldots+\lambda_\kappa\leq \kappa(\eps+\sqrt{\eps}),
$$
so it remains to find $B=B(Q)$. 

Let us consider the matrix $\tilde{Q} = \Lambda_m^{-1/2}Q\Lambda_m^{-1/2}$. Since $Q\in B_\eps(\Lambda_m)$ and $\Lambda_m$ is diagonal with all elements greater or equal than $\sqrt{\eps}$, we have $\|\tilde{Q}-I\|_\infty<  \sqrt{\eps}.$ By Gershgorin's theorem, all eigenvalues of $\tilde{Q}$ are within $m\sqrt{\eps}$ from $1$. In particular, it is invertible and we can define
\begin{equation}
B=B(Q)= \Lambda_m^{1/2}\tilde{Q}^{-1/2} \Lambda_m^{-1/2}.
\end{equation}
Using that ${Q} = \Lambda_m^{1/2}\tilde{Q} \Lambda_m^{1/2}$, it is easy to check that $B Q B^T= \Lambda_m$ and
$$
(B-I) \,Q\,(B-I)^T = \Lambda_m^{1/2}(I-\tilde{Q}^{1/2})^2 \Lambda_m^{1/2}.
$$
Since the eigenvalues of $\tilde{Q}$ are within $L\sqrt{\eps}$ from $1$, eigenvalues of $\tilde{Q}^{1/2}$ are also within $L\sqrt{\eps}$ from $1$ and, therefore, $\|I-\tilde{Q}^{1/2}\|_\infty \leq L\sqrt{\eps}.$ This implies that
$$
\tr\bigl((B-I) \,Q\,(B-I)^T \bigr) \leq L \|I-\tilde{Q}^{1/2}\|_\infty^2 \leq L\eps.
$$
Finally, since the elements of $\Lambda_m^{-1/2}$ are bounded by $\eps^{-1/4}$,
$$
\|B(Q_1)-B(Q_2)\|_\infty
\leq
L\eps^{-1/4}
\|\tilde{Q}_1^{-1/2}-\tilde{Q}_2^{-1/2}\|_\infty.
$$
Since the eigenvalues of $\tilde{Q}_1$ and $\tilde{Q}_2$ are within $L\sqrt{\eps}$ from $1$, we can take a circle of radius $1/2$ around $1$ on the complex plane, $C=\{z\in \mathbb{C} \mid |z-1|=1/2\},$ and represent
$$
\tilde{Q}_1^{-1/2}-\tilde{Q}_2^{-1/2}
=
\frac{1}{2\pi i} \int_C z^{-1/2}(z-\tilde{Q}_1)^{-1}(\tilde{Q}_2-\tilde{Q}_1)(z-\tilde{Q}_2)^{-1}\,dz,
$$
which implies that 
$$
\|\tilde{Q}_1^{-1/2}-\tilde{Q}_2^{-1/2}\|_\infty
\leq
L \|\tilde{Q}_2- \tilde{Q}_1\|_\infty
\leq
L\eps^{-1/2}  \|{Q}_2- {Q}_1\|_\infty.
$$
Combining the inequalities yields (\ref{VECdisctwotwo}) and finishes the proof.
\qed

\section{Perturbation and its consequences} \label{VECSec5label}

We will now define a direct analogue of the perturbation in the setting of the Potts spin glass \cite{PPSG} that will force the overlaps to satisfy the matrix version of the Ghirlanda-Guerra identities and all their consequences. We will first define the perturbation formally for any spin configurations, but will use it later for modifications of spin configurations defined in the previous section. For $p\geq 1$, we will use the following notation,
$$
e = (i_1,\ldots,i_p)\in \{1,\ldots, N\}^p,\,\,
\sigma_e = (\sigma_{i_1},\ldots,\sigma_{i_p})
$$
for a given $\sigma\in\Omega^N.$ Given $\lambda\in\Reals^\kappa$, we denote
$$
S_\lambda(\sigma_e) = \sum_{k\leq \kappa} \lambda_k \sigma_{i_1}(k)\cdots \sigma_{i_p}(k)
$$
and, given $n\geq 0$ and $I=(e_1,\ldots, e_{n})\in (\{1,\ldots, N\}^p)^{n}$, we let
$$
S_\lambda(\sigma_I) = S_\lambda(\sigma_{e_1})\cdots S_\lambda(\sigma_{e_n}).
$$
For integer $m\geq 1$ and $n_1,\ldots,n_m\geq 1$, let $I_j = (e_1,\ldots, e_{n_j})\in (\{1,\ldots, N\}^p)^{n_j}$ and $\lambda^j\in\Reals^\kappa$ for $1\leq j\leq m$ and consider the Hamiltonian
\begin{equation}
h_{\theta}(\sigma) = \frac{1}{N^{p(n_1+\ldots+n_m)/2}}
\sum_{I_1,\ldots,I_m} g_{I_1,\ldots,I_m}S_{\lambda^1}(\sigma_{I_1})\cdots S_{\lambda^m}(\sigma_{I_m}),
\label{VEChtheta}
\end{equation}
where $g_{I_1,\ldots,I_m}$ are standard Gaussian random variables independent for different choices of the indices. We denote the list of all parameters of the Hamiltonian by
\begin{equation}
\theta = (p,m,n_1,\ldots,n_m,\lambda^1,\ldots,\lambda^m).
\label{VECthetapar}
\end{equation}
If we recall the notation for the matrix of overlaps $R_{\ell,\ell'}$ in (\ref{VECRll}) then a straightforward calculation as in \cite{PPSG} shows that the covariance of the above Hamiltonian is given by
\begin{equation}
C^\theta_{\ell,\ell'}
=
\Cov\bigl(h_{\theta}(\sigma^\ell),h_{\theta}(\sigma^{\ell'})\bigr) 
=
\prod_{j\leq m}\bigl( R_{\ell,\ell'}^{\circ p} \lambda^j,\lambda^j \bigr)^{n_j}
\label{VECcovhtheta}
\end{equation}
for any configurations of spins $\sigma^\ell,\sigma^{\ell'}.$ Since we assume that the spins are bounded, $|\sigma_i(k)|\leq c$, the overlaps will be bounded by $c^2$ and, for $\lambda\in [-1,1]^\kappa$, we can control the quadratic form above by
$|( R_{\ell,\ell'}^{\circ p} \lambda^j,\lambda^j )| \leq \kappa^2 c^{2p}.$ 
If we denote $b_p= \kappa c^p$ then 
\begin{equation}
\bigl|C^\theta_{\ell,\ell'} \bigr|\leq b_p^{2(n_1+\ldots+n_m)}.
\label{VECCthetabound}
\end{equation}
As in \cite{PPSG}, let $\Theta$ be a collection of all $\theta$ of the type (\ref{VECthetapar}) with $p\geq 1$, $m\geq 1$, $n_1,\ldots,n_m\geq 1$, and $\lambda^1,\ldots,\lambda^m$ taking values in $([-1,1]\cap \mathbb{Q})^\kappa$ with all rational coordinates. Let us consider a one-to-one function $j_0:([-1,1]\cap \mathbb{Q})^\kappa \to\Natural$ and let
$$
j(\theta) = p+n_1+\ldots+n_m+j_0(\lambda_1)+\ldots+j_0(\lambda_m)+22m.
$$
Let $(u_{\theta})_{\theta \in \Theta}$ be i.i.d. random variables uniform on the interval $[1,2]$ and define a Hamiltonian
\begin{equation}
h_{N}(\sigma) = \sum_{\theta\in\Theta} 2^{-j(\theta)} b_p^{-(n_1+\ldots+n_m)} u_{\theta}\, h_{\theta}(\sigma).
\label{VEChNw}
\end{equation}
Conditionally on $u=(u_{\theta})_{\theta\in \Theta}$, this is a Gaussian process with the covariance
\begin{equation}
\Cov\bigl(h_{N}(\sigma^\ell),h_{N}(\sigma^{\ell'}) \bigr) 
=
 \sum_{\theta\in\Theta} 2^{-2j(\theta)} b_p^{-2(n_1+\ldots+n_m)} u_{\theta}^2
\prod_{j\leq m}\bigl( R_{\ell,\ell'}^{\circ p} \lambda^j,\lambda^j \bigr)^{n_j}.
\label{VECcovpert}
\end{equation}
In particular, the bound in (\ref{VECCthetabound}) and our choice of $j(\theta)$ imply that the variance is bounded by $1$. 

From now on, for each spin configuration $\sigma\in\Sigma^N_\eps(D)$, let $\tsigma$ denote the modified configuration $(A_\sigma \sigma_i)_{i\leq N}$ with the matrix $A_\sigma = A(R(\sigma,\sigma))$ constructed in Lemma \ref{VECLemModifyA}. Let us fix any $1/4<\gamma<1/2$, consider the sequence $s_N=N^{\gamma}$, and redefine the partition function in (\ref{VECZpartE}) by
\begin{equation}
Z_{N}(\eps,D) = \int_{\Sigma^N_\eps(D)} \exp\bigl( H_N(\sigma) + s_N h_N(\tsigma)\bigr)\, d\mu^{\otimes N}(\sigma),
\end{equation}
adding to the Hamiltonian the perturbation term $s_N h_N(\tsigma)$ depending on modified configurations.
Because the variance of $h_N$ is of order one and $\lim_{N\to\infty} N^{-1} s_N^2 = 0$, the free energy
$$
F_N\bigl(\Sigma^N_\eps(D)\bigr) = \frac{1}{N}\e\log Z_N(\eps, D)
$$
will not be affected by this perturbation in the limit. Notice that the expectation now also includes the average with respect to the uniform random variables $(u_\theta).$ One can now repeat the Aizenman-Sims-Starr calculation that leads to the representation (\ref{VECFMNAN}) with the right hand side that can be rewritten as in (\ref{VECAS2repr}), with the following minor modifications. 

First of all, the Hamiltonian $H_N'(\sigma)$ in (\ref{VECcommonH}) will be replaced by the perturbed Hamiltonian 
\begin{equation}
H_N^{\mathrm{pert}}(\sigma) = H_N'(\sigma) + s_N h_N(\tsigma),
\label{VECHpert}
\end{equation} 
and the Gibbs measure $G_N$ on $\Sigma^N_\eps(D)$ in (\ref{VECMeasureGNprime}) will be redefined by
\begin{equation}
dG_N(\sigma) = \frac{\exp H_N^{\mathrm{pert}}(\sigma)\, d\mu^{\otimes N}(\sigma)}{Z_N^{\mathrm{pert}}(\eps, D)},
\,\mbox{ where }\,
Z_N^{\mathrm{pert}}(\eps, D) = \int_{\Sigma^{N}_\eps(D)} \exp H_N^{\mathrm{pert}}(\sigma)\, d\mu^{\otimes N}(\sigma).
\label{VECMeasureGNpert}
\end{equation}
However, in the middle of this calculation the first term on the right hand side of (\ref{VECFMNAN}) will include the perturbation term $s_{N+M}h_{N+M}(\tilde{\rho})$ with $\tilde{\rho}=(A_\rho \rho_i)_{i\leq N+M}$ with the matrix $A_\rho = A(R(\rho,\rho))$ constructed in Lemma \ref{VECLemModifyA}. At that point one would like to replace it by $s_N h_N(\tsigma)$ via the interpolation
$$
\sqrt{t}s_{N+M}h_{N+M}(\tilde{\rho}) + \sqrt{1-t}s_N h_N(\tsigma)
$$
for $t\in [0,1],$ and one needs to check that this introduces an error that vanishes as $N\to\infty.$ If, conditionally on $(u_\theta)$, we think of the right hand side of (\ref{VECcovpert}) as a function of the overlap matrix $R_{\ell,\ell'}$, denote it by $f(R_{\ell,\ell'})$ and compute the derivative of the first term on the right hand side of (\ref{VECFMNAN}) in the parameter $t$ in the above interpolation using Gaussian integration by parts, we will see that the order of the derivative will be determined by the quantities of the type
$$
(N+M)^{2\gamma}f\bigl(R(\tilde{\rho}^1,\tilde{\rho}^2)\bigr)-N^{2\gamma}f\bigl(R(\tsigma^1,\tsigma^2)\bigr)
$$
(see e.g. Section 3.5 in \cite{SKmodel} for details). Let us recall that we write the configuration $\rho\in \Omega^{N+M}$ as $(\sigma,\tau)$ for $\sigma\in\Omega^N$ and $\tau\in\Omega^M,$ and
$$
R(\rho,\rho) = \frac{N}{N+M}R(\sigma,\sigma)+\frac{M}{N+M}R(\tau,\tau).
$$
For a fixed $M$, this implies that $|R(\rho,\rho) - R(\sigma,\sigma)|= \bigO{N^{-1}}$ so, from the equation (\ref{VECRoverR}) and Lemma \ref{VECLemModifyA}, we see that $|R(\tilde{\rho}^1,\tilde{\rho}^2) - R(\tsigma^1,\tsigma^2)|= \bigO{(N\eps)^{-1}}$. Since $(N+M)^{2\gamma} - N^{2\gamma}$ is of the order $N^{-(1-2\gamma)}$ and the derivative of $f$ is bounded on compacts uniformly over $(u_\theta)$, the order of the derivative in the above interpolation will be $N^{-(1-2\gamma)}/\eps$ and the error introduced by the interpolation will vanish in the limit $N\to\infty.$

As in the Potts spin glass model in \cite{PPSG}, the perturbation term $s_N h_N(\tsigma)$ is introduced to ensure the validity of some overlap-matrix version of the classical Ghirlanda-Guerra identities \cite{GG} for the Gibbs measure (\ref{VECMeasureGNpert}). Given replicas $(\sigma^\ell)$ from the Gibbs measure $G_N$ on $\Sigma_\eps^N(D),$ let us denote by 
\begin{equation}
\tilde{R}_{\ell,\ell'}=R(\tsigma^{\ell},\tsigma^{\ell'})
\,\mbox{ and }\,
\tilde{R}^n = \bigl(\tilde{R}_{\ell, \ell'}\bigr)_{\ell,\ell'\leq n}
\label{VECtildeOR}
\end{equation}
for any $n\geq 2$. Similarly to (\ref{VECcovhtheta}), let us define
\begin{equation}
\tilde{C}^\theta_{\ell,\ell'}
=
\Cov(h_{\theta}(\tsigma^\ell),h_{\theta}(\tsigma^{\ell'})) 
=
\prod_{j\leq m}\bigl( \tilde{R}_{\ell,\ell'}^{\circ p} \lambda^j,\lambda^j \bigr)^{n_j}.
\label{VECtcovhtheta}
\end{equation}
Consider an arbitrary bounded measurable function $f=f(\tilde{R}^n)$ and, for $\theta \in \Theta$, let 
\begin{equation}
\varDelta(f,n,\theta) = 
\Bigl|
{\e}  \bigl\la f \tilde{C}_{1,n+1}^{\theta} \bigr\ra -  \frac{1}{n} {\e} \bigl\la f \bigr\ra {\e} \bigl\la \tilde{C}_{1,2}^{\theta} \bigr\ra - \frac{1}{n}\sum_{\ell=2}^{n}{\e} \bigl\la f \tilde{C}_{1,\ell}^{\theta} \bigr\ra
\Bigr|,
\label{VECGGfinite}
\end{equation}
where ${\e}$ denotes the expectation conditionally on the i.i.d. uniform sequence $u=(u_{\theta})_{\theta\in \Theta}$. If we denote by $\e_u$ the expectation with respect to $u$ then the following holds. 
\begin{lemma}\label{VECThGG} 
For any $n\geq 2$ and any bounded measurable function $f=f(\tilde{R}^n)$, for all $\theta\in\Theta$,
\begin{equation}
\lim_{N\to\infty} \e_u \smsp \varDelta(f,n,\theta) = 0.
\label{VECGGxlim}
\end{equation} 
\end{lemma}
\textbf{Proof.} The proof is identical to proof of Theorem 3.2 in \cite{SKmodel}, but we should emphasize one more time why we defined the perturbation Hamiltonian in terms of modified configurations $\tsigma.$ The reason is because the proof of the equation (\ref{VECGGxlim}) follows from some Gaussian integration by parts computation involving one term $h_{\theta}(\tsigma)$ in the perturbation (\ref{VEChNw}), but this computation only works if the covariance $\tilde{C}_{\ell,\ell}^{\theta}$ corresponding to the same configuration $\tsigma^\ell$ is constant independent of the configuration. Otherwise, some additional terms will appear. By the construction of the modified configurations in Lemma \ref{VECLemModifyA}, 
$$
\tilde{C}_{\ell,\ell}^{\theta} = 
\prod_{j\leq m}\bigl( \tilde{R}_{\ell,\ell}^{\circ p} \lambda^j,\lambda^j \bigr)^{n_j}
=
\prod_{j\leq m}\Bigl( D_\eps^{\circ  p} \lambda^j,\lambda^j \Bigr)^{n_j}
$$
are, indeed, independent of the configuration. Without spin modification, the self-overlap $R(\sigma^\ell,\sigma^\ell)$ would be non-constant, since it is only constrained to be in the $\eps$-neighbourhood of $D\in\DD.$ With the small modification of spins that fixes the self-overlap to be equal to $D_\eps$, the proof of the Ghirlanda-Guerra identities goes through without any changes.
\qed

\medskip
Let us now summarize main consequences of this result obtained in \cite{PPSG}. Using (\ref{VECGGxlim}), one can choose a non-random sequence $u^N=(u^N_{\theta})_{\theta\in \Theta}\in [1,2]^\Theta$ such that 
\begin{equation}
\lim_{N\to\infty} \smsp \varDelta(f,n,\theta) = 0
\,\mbox{ for all }\, \theta\in\Theta
\label{VECGGxlim2}
\end{equation} 
for the Gibbs measure $G_N$ with the parameters $u$ in the perturbation (\ref{VEChNw}) equal to $u^N$ rather than random. Consider any such sequence $u^N$ and consider any subsequence $(N_k)_{k\geq 1}$ along which the array $(\tilde{R}_{\ell,\ell'})_{\ell,\ell'\geq 1}$ of overlap matrices in (\ref{VECtildeOR}) converges in distribution under the measure $\e G_N^{\otimes\infty}$. Let us slightly abuse notation and denote the limiting array again by
\begin{equation}
R_{\ell,\ell'} = \bigl(R_{\ell,\ell'}^{k,k'} \bigr)_{k,k'\leq \kappa},\,
R^n = \bigl(R_{\ell, \ell'}\bigr)_{\ell,\ell'\leq n},\,
C^\theta_{\ell,\ell'} =
\prod_{j\leq m} \bigl( R_{\ell,\ell'}^{\circ p} \lambda^j,\lambda^j \bigr)^{n_j}.
\label{VECRwlim}
\end{equation}
Then the equations (\ref{VECGGfinite}) and (\ref{VECGGxlim2}) imply that \begin{equation}
\e f(R^n) C^\theta_{1,n+1}
=
 \frac{1}{n}\e f(R^n)  \e C^\theta_{1,2}
+ \frac{1}{n}\sum_{\ell=2}^{n}\e f(R^n) C^\theta_{1,\ell}
\label{VECGGwp}
\end{equation}
for all $\theta\in\Theta$. Since $C^\theta_{\ell,\ell'}$ is a continuous function of $\lambda^j\in [-1,1]^{\kappa}$ for $j\leq m$, (\ref{VECGGwp}) holds a posteriori for all values of $\lambda^j$, not only with rational coordinates. 

For any $p\geq 1$, $\lambda^1,\ldots,\lambda^m \in [-1,1]^{\kappa}$ and a bounded measurable function $\varphi\colon \Reals^m\to\Reals$, let
\begin{equation}
Q_{\ell,\ell'} = \varphi \Bigl( \bigl( R_{\ell,\ell'}^{\circ p} \lambda^1,\lambda^1 \bigr),\ldots, \bigl( R_{\ell,\ell'}^{\circ p} \lambda^m,\lambda^m \bigr) \Bigr).
\label{VECQphi}
\end{equation} 
As in Theorem 2 in \cite{PPSG}, the next result immediately follows from (\ref{VECGGwp}).
\begin{theorem}\label{VECThGGms}
For any $n\geq 2$ and any bounded measurable function $f=f(R^n)$,
\begin{equation}
\e f(R^n) Q_{1,n+1}
=
 \frac{1}{n}\e f(R^n)  \smsp \e Q_{1,2}
+ \frac{1}{n}\sum_{\ell=2}^{n}\e f(R^n) Q_{1,\ell}.
\label{VECGGms}
\end{equation}
\end{theorem}
In addition to well-known standard consequences of the classical Ghirlanda-Guerra identities, which are contained in (\ref{VECGGms}), the main consequence about the structure of the limiting overlap arrays was the following result in Theorem 3 in \cite{PPSG} about the synchronization of the blocks of overlaps.

\begin{theorem}\label{VECTh2}
If the overlap array satisfies (\ref{VECGGms}) for all choices of parameters then
\begin{equation}
R_{\ell,\ell'} = \Phi\bigl(\tr(R_{\ell,\ell'}) \bigr) \mbox{ a.s.}
\end{equation}
for some function $\Phi\colon \Reals^+ \to \Gamma_\kappa$, which is non-decreasing in $\Gamma_\kappa$,  $\Phi(x')-\Phi(x)\in \Gamma_\kappa$ for all $x\leq x',$ and Lipschitz continuous, $\|\Phi(x') - \Phi(x)\|_1 \leq L |x'-x|$.
\end{theorem}
We now return to the cavity computation and explain the next steps.

\section{Cavity computation, part 2}\label{VECSec6label}

If we denote the quantity in (\ref{VECAS2repr}) by
\begin{align}
A_{N,M} = &\,\,
\frac{1}{M}
\e \log \Bigl\la
\int_{\Sigma^{M}_\eps(D)} \exp\Bigl( \sum_{i\leq M}\sum_{k\leq\kappa} \tau_i(k) Z_{i}^k(\sigma)
\Bigr)\, d\mu^{\otimes M}(\tau)
\Bigr\ra_N
\nonumber
\\
&
-
\frac{1}{M}
\e \log \Bigl\la \exp \sqrt{M} Y(\sigma) \Bigr\ra_N
\label{VECASreprpart2}
\end{align}
then in the previous section we explained that the lower bound
\begin{align}
\liminf_{N\to\infty} F_N\bigl(\Sigma^N_\eps(D) \bigr)
\geq
\liminf_{N\to\infty} A_{N,M}
\label{VECFMNANpart2}
\end{align}
holds for the Gibbs measure $G_N$ in (\ref{VECMeasureGNpert}) corresponding to the perturbed Hamiltonian.
Recall that in this case the expectation $\e$ in (\ref{VECASreprpart2}) includes the average $\e_u$ in the uniform random variables $u=(u_{\theta})_{\theta\in\Theta}$ in the definition of the perturbation Hamiltonian (\ref{VEChNw}). By Lemma 3.3 in \cite{SKmodel}, one can choose a non-random sequence $u^N=(u^N_{\theta})_{\theta\in \Theta}$ changing with $N$ such that both (\ref{VECGGxlim2}) and (\ref{VECFMNANpart2}) hold for the Gibbs measure $G_N$ with the parameters $u$ in the perturbation Hamiltonian (\ref{VEChNw}) equal to $u^N$ rather than random. 

Next, similarly to (\ref{VECASreprpart2}), let us define
\begin{align}
\tilde{A}_{N,M} = &\,\,
\frac{1}{M}
\e \log \Bigl\la
\int_{\Sigma^{M}_\eps(D)} \exp\Bigl( \sum_{i\leq M}\sum_{k\leq\kappa} \tau_i(k) Z_{i}^k(\tsigma)
\Bigr)\, d\mu^{\otimes M}(\tau)
\Bigr\ra_N
\nonumber
\\
&
-
\frac{1}{M}
\e \log \Bigl\la \exp \sqrt{M} Y(\tsigma) \Bigr\ra_N,
\label{VECtASreprpart2}
\end{align}
where we replaced the configuration $\sigma$ which indexes the processes $Z_i^k$ and $Y$ by the modified configuration $\tsigma$ defined in Lemma \ref{VECLemModifyA}. As in (\ref{VECCovz}), (\ref{VECCovy}) and (\ref{VECCovysum}), up to smaller order terms which we can omit, the covariance of these processes indexed by modified configurations is given by
\begin{align}
\e Z_i^k(\tsigma^{\ell}) Z_i^{k'}(\tsigma^{\ell'})
&= 
\xi_{k,k'}'\bigl(\tilde{R}_{\ell,\ell'}^{k,k'} \bigr),
\label{VECCovzpart2}
\\
\e Y(\tsigma^\ell) Y(\tsigma^{\ell'}) 
&= 
\su\bigl(\theta(\tilde{R}_{\ell,\ell'})\bigr).
\label{VECCovypart2}
\end{align}
By (\ref{VECoverchange}) and (\ref{VECdiscone}), $\|\tilde{R}_{\ell,\ell'}-{R}_{\ell,\ell'}\|_\infty\leq L\eps^{1/4}$ so the covariance of these processes is affected only slightly by this substitution. In particular, using the usual Gaussian interpolation of the form 
$$
\sqrt{t}Z_{i}^k(\sigma) + \sqrt{1-t} Z_{i}^k(\tsigma),\,\,
\sqrt{t}Y^k(\sigma) + \sqrt{1-t} Y^k(\tsigma),
$$
one can show that $|\tilde{A}_{N,M}-{A}_{N,M}|\leq L\eps^{1/4}$ and, therefore,
\begin{align}
\liminf_{N\to\infty} F_N\bigl(\Sigma^N_\eps(D) \bigr)
\geq
\liminf_{N\to\infty} \tilde{A}_{N,M} - L\eps^{1/4}.
\label{VECFMNANpart22}
\end{align}
Let us take a subsequence along which the lower limit on the right hand side is achieved and then pass to another subsequence along which the distribution of the array $(\tilde{R}_{\ell,\ell'})_{\ell,\ell'\geq 1}$ under $\e G_N^{\otimes\infty}$ converges. Let us denote the array with the limiting distribution by $(R^{M}_{\ell,\ell'})_{\ell,\ell'\geq 1}$, because the limit was taken for a fixed $M$ and may depend on $M$. Notice that, because of the definition of the modified configurations $\tsigma$, the diagonal overlap blocks are fixed,
\begin{equation}
\tilde{R}_{\ell,\ell'} = \tilde{R}^M_{\ell,\ell'} = D_\eps.
\label{VECconstraintDeps}
\end{equation}
As in the case of the Potts spin glass in \cite{PPSG}, we now recall the well-known fact (see e.g. the proof of Theorem 1.3 in \cite{SKmodel}) that both terms in (\ref{VECtASreprpart2}) are continuous functionals of the distribution of the array $(\tilde{R}_{\ell,\ell'})_{\ell,\ell'\geq 1}$ under $\e G_N^{\otimes\infty}$, so to describe the limit we need to understand how this functional looks like for the limiting array  $(R^{M}_{\ell,\ell'})_{\ell,\ell'\geq 1}$. We showed that, due to the perturbation of the Hamiltonian, this array satisfies the generalized Ghirlanda-Guerra identities in Theorem \ref{VECThGGms} and the synchronization property in Theorem \ref{VECTh2}. Moreover, by Theorem \ref{VECThGGms}, the array $(\tr(R^{M}_{\ell,\ell'}))_{\ell,\ell'\geq 1}$ itself satisfies the classical Ghirlanda-Guerra identities and, by the results in Chapter 2 of \cite{SKmodel}, it can be generated by the Ruelle probability cascades. This means that the proof can be finished exactly as in \cite{PPSG} if we can only show the Lipschitz continuity and decoupling properties of the analogues of the functionals in (\ref{VECtASreprpart2}) for the Ruelle probability cascades, which we will  do next.

\section{Functionals of the Ruelle probability cascades}\label{VECSec7label}

Let us consider a discrete path $\pi\in \Pi_{\Delta}$ defined as in (\ref{VECgammapath}) in terms of the sequences
\begin{equation}
x_{-1}=0< x_0 < \ldots < x_{r-1} < x_r = 1
\label{VECxsRPC}
\end{equation}
and a monotone sequence of Gram matrices in $\Gamma_\kappa$,
\begin{equation}
0=\gamma_0\leq \gamma_1 \leq \ldots \leq \gamma_{r-1}\leq \gamma_r = \Delta,
\label{VECgammasRPC}
\end{equation}
only now the final constraint is given by some arbitrary $\Delta\in \Gamma_\kappa$. Let us consider the Gaussian processes $Z(\alpha)$ and $Y(\alpha)$ defined as in Section \ref{VECSec2label} with the covariances
\begin{align}
\Cov\bigl(Z(\alpha^1), Z(\alpha^2)\bigr) &= \xi'(\gamma_{\alpha^1\wedge\alpha^2}),
\nonumber
\\
\Cov\bigl(Y(\alpha^1), Y(\alpha^2)\bigr) &= \su(\theta(\gamma_{\alpha^1\wedge\alpha^2})),
\label{VECCDRPC}
\end{align}
and let $Z_i(\alpha)$ be independent copies of $Z(\alpha)$ for $i\geq 1.$ The path $\pi$, including the constraint $\Delta$ in (\ref{VECgammasRPC}), will be fixed for the rest of this section so we will not write the dependence on it explicitly. Let us consider the analogues of the functionals in (\ref{VECASreprpart2}), 
\begin{align}
f_M^1\bigl(B_\eps(D) \bigr)
&= 
\frac{1}{M}
\e \log \sum_{\alpha\in\Natural^r} v_\alpha 
\int_{\Sigma^M_\eps(D)} \exp\Bigl( \sum_{i\leq M}\sum_{k\leq\kappa} \tau_i(k) Z_{i}^k(\alpha)
\Bigr)\, d\mu^{\otimes M}(\tau),
\label{VECfN1}
\\
f_M^2
&= 
\frac{1}{M}
\e \log \sum_{\alpha\in\Natural^r} v_\alpha  \exp \sqrt{M} Y(\alpha).
\label{VECfN2}
\end{align}
Later we will replace the final constraints $\Delta$ in (\ref{VECgammasRPC}) by $D_\eps$ defined in (\ref{VECQLQeps}), but in this section we will let $\Delta$ be arbitrary and unrelated to the constraint on the configurations $\tau\in\Sigma^M_\eps(D).$ The functionals (\ref{VECfN1}) and (\ref{VECfN2}) are precisely the functionals that appeared at the end of Guerra's replica symmetry breaking interpolation in Section \ref{VECSec2label} (only now we write $M$ instead of $N$, $\tau$ instead of $\sigma$, and $\Delta$ instead of $D$ in (\ref{VECgammasRPC})). We have seen in the proof of Lemma \ref{VECLemAlmostUp} that
\begin{equation}
f_M^2
=
\frac{1}{2}
\sum_{0\leq j\leq r-1} x_j  \su\bigl( \theta(\gamma_{j+1}) - \theta(\gamma_{j})\bigr).
\label{VECfunp2RPC}
\end{equation}
If we recall the functional $\Phi(\lambda)=\Phi(\lambda, \Delta, \pi)$ defined in (\ref{VECPhi}) (with $D$ now replaced by $\Delta$), in the proof of Lemma \ref{VECLemAlmostUp} we appealed to the properties of the Ruelle probability cascades to claim that
\begin{equation}
\Phi(\lambda) 
=
\e \log \sum_{\alpha\in\Natural^r} v_\alpha \int_{\Omega} \exp
\Bigl( \sum_{k\leq \kappa} \tau_1(k) Z_{i}^k(\alpha) + \sum_{k\leq k'} \lambda_{k,k'} \tau_1(k)\tau_1(k') \Bigr)\, d\mu(\tau_1).
\label{VECPhilambda}
\end{equation}
We also showed there that, for any $\lambda = (\lambda_{k,k'})_{1\leq k\leq k'\leq \kappa}\in \Reals^{\kappa(\kappa+1)/2}$,
$$
f_M^1\bigl(B_\eps(D) \bigr)
\leq \eps\|\lambda\|_1 -\sum_{k\leq k'} \lambda_{k,k'} D_{k,k'} + \Phi(\lambda).
$$
We will now show that, if we omit $\eps\|\lambda\|_1$ and minimize over $\lambda$, this bound becomes asymptotically sharp. For $D\in\DD,$ let us denote
\begin{equation}
\Phi^*(D): = \inf_{\lambda}\Bigl(  -\sum_{k\leq k'} \lambda_{k,k'} D_{k,k'} + \Phi(\lambda) \Bigr).
\label{VECPhistar}
\end{equation}
Next lemma will follow by an adaptation of a standard smoothing technique (see e.g. Section 2.2.2 in \cite{Dembo}), combined with some straightforward spin glass calculations. 
\begin{lemma}\label{VECLem2}
For any $\eps>0$ and $D\in\DD,$
\begin{equation}
\liminf_{M\to\infty} f_M^1\bigl(B_\eps(D) \bigr) \geq \Phi^*(D).
\label{VECfNlim}
\end{equation}
\end{lemma}

Let $g_i=(g_i(k,k'))_{k,k'\leq \kappa}$ be i.i.d. symmetric $\kappa\times\kappa$ matrices with independent Gaussian entries with variance $\delta>0$ except for the symmetry constraint $g_i(k,k')=g_i(k',k)$, and let $g=(g_1,\ldots,g_M).$ We will denote the distribution of $g_i$ on $\Reals^{\kappa\times\kappa}$ by $\nu$. Define $\xoverline{g}=M^{-1}\sum_{i\leq M} g_i$ and, for any subset $S\subseteq \Reals^{\kappa\times\kappa}$, let us consider the set
\begin{equation}
\Sigma(S) = \Bigl\{ (\tau,g) \mid R(\tau,\tau) + \xoverline{g} \in S \Bigr\}.
\end{equation}
Similarly to (\ref{VECfN1}), let us define
\begin{equation}
f_M^g(S)
= 
\frac{1}{M}
\e \log \sum_{\alpha\in\Natural^r} v_\alpha 
\int_{\Sigma(S)} \exp\Bigl( \sum_{i\leq M}\sum_{k\leq\kappa} \tau_i(k) Z_{i}^k(\alpha)
\Bigr)\, d\mu^{\otimes M}(\tau) d\nu^{\otimes M}(g).
\label{VECfN12}
\end{equation}
Without the Gaussian random variables $g$ and with $S=B_\eps(D)$, this would be exactly $f_M^1(B_\eps(D)).$ Similarly to (\ref{VECPhilambda}), let us define
\begin{align}
\Phi_g(\lambda) 
=
\e \log \sum_{\alpha\in\Natural^r} v_\alpha \int \exp
&
\Bigl( \,
\sum_{k\leq \kappa} \tau_1(k) Z_{i}^k(\alpha) 
\nonumber
\\
&\,\,\,
+ \sum_{k\leq k'} \lambda_{k,k'} \bigl(\tau_1(k)\tau_1(k') +g_1(k,k')\bigr)\Bigr)\, d\mu(\tau_1)d\nu(g_1).
\label{VECPhiGlambda}
\end{align}
Since $g_1(k,k')$ for $k\leq k'$ are independent Gaussian with variance $\delta,$
\begin{equation}
\Phi_g(\lambda) 
=
\Phi(\lambda) + \frac{\delta}{2}\sum_{k\leq k'}\lambda_{k,k'}^2.
\label{VECPhiGlambdaD}
\end{equation}
Next, as in (\ref{VECPhistar}), let us define
\begin{equation}
\Phi_g^*(D): = \inf_{\lambda}\Bigl(  -\sum_{k\leq k'} \lambda_{k,k'} D_{k,k'} + \Phi_g(\lambda) \Bigr).
\label{VECPhiGstar}
\end{equation}
Since the symmetric random matrix $\xoverline{g}$ is not necessarily positive-semidefinite, let us redefine the set $B_\eps(D)$ in (\ref{VECBepsD}) to be a subset
\begin{equation}
B_\eps(D) = \Bigl\{ \gamma\in \Reals^{\kappa\times\kappa} \mid \|\gamma-D\|_\infty <\eps \Bigr\}
\label{VECBepsD2}
\end{equation}
of $\Reals^{\kappa\times\kappa}$ rather than $\Gamma_\kappa$. We will begin by proving the following.
\begin{lemma}\label{VECLem2G}
For any $\eps>0$ and $D\in\DD,$
\begin{equation}
\liminf_{M\to\infty} f_M^g\bigl(B_\eps(D)\bigr) \geq \Phi_g^*(D).
\label{VECfNlim2}
\end{equation}
\end{lemma}
\textbf{Proof.} Since $\Phi(\lambda)$ is convex and grows at most linearly in $\lambda$, the presence of the quadratic second term in (\ref{VECPhiGlambdaD}) guarantees that the infimum in (\ref{VECPhiGstar}) is achieved on some critical point $\lambda$ such that 
\begin{equation}
\nabla\Phi_g(\lambda) = D.
\label{VECnablaPhi}
\end{equation}
Here $D$ represent only the upper half of the symmetric matrix $D$, but we will abuse the notation and simply write $D.$ In other words, with this choice of $\lambda$,
$$
\Phi_g^*(D): = -\sum_{k\leq k'} \lambda_{k,k'} D_{k,k'} + \Phi_g(\lambda).
$$
Similarly to (\ref{VECfN12}), consider the functional
\begin{align}
f_M^g(S,\lambda)
= 
\frac{1}{M}
\e \log \sum_{\alpha\in\Natural^r} v_\alpha 
\int_{\Sigma(S)} \exp
&
\Bigl( \,
\sum_{i\leq M}\sum_{k\leq\kappa} \tau_i(k) Z_{i}^k(\alpha)
\nonumber
\\
&\,\,
+ M \sum_{k\leq k'} \lambda_{k,k'}(R(\tau,\tau)+\xoverline{g})_{k,k'}
\Bigr)\, d\mu^{\otimes M}(\tau) d\nu^{\otimes M}(g).
\label{VECfN13}
\end{align}
Notice that, by the standard properties of the Ruelle probability cascades that were already invoked in the proof of the Guerra upper bound,
\begin{equation}
f_M^g(\Reals^{\kappa\times\kappa},\lambda)=\Phi_g(\lambda)
\label{VECtriuppe}
\end{equation}
with $\Phi_g$ defined in (\ref{VECPhiGlambda}). Let us now consider the complement of $B_\eps(D)$ in $\Reals^{\kappa\times \kappa}$ and let us cover it by half-spaces of the form
$$
H_{k,k'}^+=\Bigl\{x\in \Reals^{\kappa\times \kappa} \mid x_{k,k'}\geq D_{k,k'} +\eps \Bigr\},\,\,
H_{k,k'}^-=\Bigl\{x\in \Reals^{\kappa\times \kappa} \mid x_{k,k'}\leq D_{k,k'} -\eps \Bigr\}.
$$
Because all the matrices are symmetric, we only need to consider indices $k\leq k'$. Let us consider one such half-space, for example, $H=H^+_{m,m'}.$ Let us denote  
$$
e_{m,m'}=\Bigl(\I\bigl((k,k')=(m,m')\bigr) \Bigr)_{k\leq k'}.
$$
Since, for $t\geq 0$ and $x\in H,$ we have $t(x_{m,m'} - D_{m,m'}-\eps)\geq 0,$ this together with (\ref{VECtriuppe}) implies that
\begin{equation}
f_M^g(H,\lambda)\leq -t(D_{m,m'}+\eps) + \Phi_g\bigl(\lambda+t e_{m,m'} \bigr).
\label{VECtupper}
\end{equation}
For $t=0$, this upper bound equals $\Phi_g(\lambda)$ and, by (\ref{VECnablaPhi}), the derivative
$$
\frac{\partial}{\partial t} \Phi_g\bigl(\lambda+t e_{m,m'} \bigr) \Bigr|_{t=0}
=
\frac{\partial}{\partial \lambda_{m,m'}} \Phi_g\bigl(\lambda \bigr) = D_{m,m'},
$$
so the derivative of the right hand side of (\ref{VECtupper}) at $t=0$ equals $-\eps.$ It is tedious but straightforward to check that the second derivatives of $\Phi_g$ are bounded on compacts, and as a result,
$$
f_M^g(H,\lambda)\leq \Phi_g\bigl(\lambda\bigr) - \eps t + \frac{Lt^2}{2}
$$
for $t\in [0,1]$ for some large enough $L.$ For $t=\eps/L$ this yields the bound
\begin{equation}
f_M^g(H,\lambda)\leq \Phi_g\bigl(\lambda\bigr) - \frac{\eps^2}{2L}.
\label{VECFMgH}
\end{equation}
The same bound can be obtained by a similar argument for any $H=H^-_{m,m'}.$ The argument in the proof of Lemma 6 in \cite{PPSG} (which we do not repeat here) shows that if $A_j(\alpha)$ for $1\leq j\leq n$ are some positive functions of the Gaussian processes $Z_i(\alpha)$ and $Y(\alpha)$ then
\begin{equation}
\e \log \sum_{\alpha\in\Natural^r} v_\alpha \sum_{j\leq n} A_j(\alpha)
\leq
\frac{\log n}{x_0} + \max_{j\leq n} \e \log \sum_{\alpha\in\Natural^r} v_\alpha A_j(\alpha),
\label{VECRPCboundA}
\end{equation}
where $x_0$ is the first element in the sequence (\ref{VECxsRPC}). Since $\Reals^{\kappa\times\kappa}$ can be covered by $B_\eps(D)$ and the half-spaces as above, this implies that
$$
\Phi_g(\lambda)
=
f_M^g\bigl(\Reals^{\kappa\times\kappa},\lambda\bigr)
\leq
\frac{\log(2\kappa+1)}{Mx_0}+ \max\Bigl(
f_M^g(B_\eps(D),\lambda),
\max_{H} f_M^g(H,\lambda)
\Bigr).
$$
The maximum $\max_H$ on the right hand side is over the above half-spaces and the bound (\ref{VECFMgH}) ensures that one can not have
$$
\Phi_g(\lambda)
\leq 
\frac{\log(2\kappa+1)}{Mx_0}+ \max_{H} f_M^g(H,\lambda),
$$
for large $M$. Therefore, we must have
$$
\Phi_g(\lambda)
\leq
\frac{\log(2\kappa+1)}{Mx_0}+ f_M^g(B_\eps(D),\lambda).
$$
On the other hand, from the definition of these functionals,
$$
 f_M^g\bigl(B_\eps(D),\lambda\bigr) \leq  \sum_{k\leq k'} \lambda_{k,k'} D_{k,k'} + f_M^g\bigl(B_\eps(D)\bigr) + \eps \|\lambda\|_1.
$$
The above two inequalities imply that
$$
\liminf_{M\to\infty}  f_M^g\bigl(B_\eps(D) \bigr) \geq -\sum_{k\leq k'} \lambda_{k,k'} D_{k,k'} + \Phi_g(\lambda) - \eps \|\lambda\|_1 = \Phi_g^*(D) - \eps \|\lambda\|_1.
$$
Since $f_M^g(B_\eps(D))$ is increasing in $\eps$, we can let $\eps\downarrow 0$ on the right hand side while fixing it on the left hand side. This finishes the proof.
\qed

\medskip
We can now deduce Lemma \ref{VECLem2} from Lemma \ref{VECLem2G}.

\medskip
\noindent
\textbf{Proof of Lemma \ref{VECLem2}.}
Using that
$$
\Bigl\{R(\tau,\tau) + \xoverline{g} \in B_\eps(D) \Bigr\}
\subseteq
\Bigl\{R(\tau,\tau) \in B_{2\eps}(D) \Bigr\}
\bigcup
\Bigl\{ \|\xoverline{g}\|_\infty \geq \eps\Bigr\},
$$
we can bound
$$
\int_{\Sigma(B_\eps(D))} \exp
\Bigl( \,
\sum_{i\leq M}\sum_{k\leq\kappa} \tau_i(k) Z_{i}^k(\alpha)
\Bigr)\, d\mu^{\otimes M}(\tau) d\nu^{\otimes M}(g)
\leq 
A_1(\alpha) + A_2(\alpha)
$$
where
\begin{align*}
A_1(\alpha)
& =
\int_{\Sigma_{2\eps}^M(D)} \exp
\Bigl( \,
\sum_{i\leq M}\sum_{k\leq\kappa} \tau_i(k) Z_{i}^k(\alpha)
\Bigr)\, d\mu^{\otimes M}(\tau),
\\
A_2(\alpha)
&= 
\p\bigl(\|\xoverline{g}\|_\infty \geq \eps\bigr)
\int_{\Omega^M} \exp
\Bigl( \,
\sum_{i\leq M}\sum_{k\leq\kappa} \tau_i(k) Z_{i}^k(\alpha)
\Bigr)\, d\mu^{\otimes M}(\tau).
\end{align*}
Using this and (\ref{VECRPCboundA}), we can bound
\begin{equation}
f_M^g\bigl(B_\eps(D)\bigr)
\leq
\frac{\log 2}{Mx_0} + \max\Bigl(f_M^1\bigl(B_{2\eps}(D)\bigr), A_2\Bigr),
\label{VECflemal}
\end{equation}
where
\begin{align*}
A_2 
&=
 \frac{1}{M}\log \p\bigl(\|\xoverline{g}\|_\infty \geq \eps\bigr)
 +
 \frac{1}{M}\e\log 
 \sum_{\alpha\in\Natural^r} v_\alpha
 \int_{\Omega^M} \exp
\Bigl( \,
\sum_{i\leq M}\sum_{k\leq\kappa} \tau_i(k) Z_{i}^k(\alpha)
\Bigr)\, d\mu^{\otimes M}(\tau).
\end{align*}
Since $\xoverline{g}$ is a vector of $\kappa$ independent Gaussian random variables with the variance $\delta/M,$ if $z$ is a standard Gaussian random variable,
$$
 \frac{1}{M}\log \p\bigl(\|\xoverline{g}\|_\infty \geq \eps\bigr)
 =
\frac{\kappa}{M}\log \p\Bigl(|z|\geq \frac{\eps\sqrt{M}}{\sqrt{\delta}} \Bigr)
\leq 
-\frac{k\eps^2}{2\delta}+\frac{\kappa\log 2}{M}.
$$
The second term in $A_2$ is bounded by some constant, which can be seen by taking the expectation inside the logarithm. By letting $\delta\downarrow 0$, one can make $A_2\to-\infty.$ On the other hand, (\ref{VECPhiGlambdaD}) implies that $\Phi_g^*(D)\geq \Phi^*(D)$ and, by the previous lemma,
$$
\liminf_{M\to\infty} f_M^g\bigl(B_\eps(D)\bigr) \geq \Phi^*(D).
$$
Therefore, letting $M\to\infty$ and then $\delta\downarrow 0$ in (\ref{VECflemal}) finishes the proof.
\qed

\medskip
In addition to the above lower bound, we need to recall standard Lipschitz continuity property for the functionals $f_M^1(B_\eps(D), \pi)$ in (\ref{VECfN1}) and $f_M^2(\pi)$ in (\ref{VECfN2}), where we now make the dependence on the path $\pi\in \Pi_\Delta$ defined in terms of the sequences (\ref{VECxsRPC}) and (\ref{VECgammasRPC}) explicit. 
\begin{lemma}\label{VECLemLC1}
For any $\Delta\in \Gamma_\kappa$ and for any two discrete paths $\pi, \tilde{\pi}\in \Pi_\Delta$,
\begin{equation}
\bigl|f_M^1(B_\eps(D), \pi)-f_M^1(B_\eps(D), \tilde{\pi})\bigr|
\leq
L \int_{0}^{1}\! \|\pi(x) - \tilde{\pi}(x)\|_1 \, dx
\label{VECLipschitz1}
\end{equation}
and 
\begin{equation}
\bigl|f_M^2(\pi) -f_M^2(\tilde{\pi})\bigr|
\leq
L \int_{0}^{1}\! \|\pi(x) - \tilde{\pi}(x)\|_1 \, dx.
\label{VECLipschitz2}
\end{equation}
\end{lemma}
The representation of $f_M^2(\pi)$ in (\ref{VECfunp2RPC}), especially when written in the form (\ref{VECrearrange}), makes the second equation (\ref{VECLipschitz2}) obvious. The proof of the first one is identical to the proof of Lemma 7 in \cite{PPSG}.

\section{Cavity computation, part 3}\label{VECSec8label}

The rest of the proof is almost identical to the proof of the lower bound in \cite{PPSG}, and we will only sketch it here without repeating all the details. We showed in (\ref{VECFMNANpart22}) that 
\begin{align}
\liminf_{N\to\infty} F_N\bigl(\Sigma^N_\eps(D)\bigr)
\geq
\liminf_{N\to\infty} \tilde{A}_{N,M} - L\eps^{1/4}
\label{VECFMNANpart22again}
\end{align}
with $\tilde{A}_{N,M}$ defined in (\ref{VECtASreprpart2}). We denoted by $(R^{M}_{\ell,\ell'})_{\ell,\ell'\geq 1}$  the limit of the array $(\tilde{R}_{\ell,\ell'})_{\ell,\ell'\geq 1}$ in (\ref{VECtildeOR}) in distribution over some subsequence of the sequence along which the lower limit in (\ref{VECFMNANpart22again}) is achieved. One can then take the lower limit of (\ref{VECFMNANpart22again}) as $M\to\infty$ and choose a subsequence along which  $(R^{M}_{\ell,\ell'})_{\ell,\ell'\geq 1}$ converges in distribution to some array  $(R^{\infty}_{\ell,\ell'})_{\ell,\ell'\geq 1}$. All these arrays satisfy the Ghirlanda-Guerra identities in Theorem \ref{VECThGGms}. Using the synchronization property in Theorem \ref{VECTh2} and well-known approximation properties of arrays satisfying the Ghirlanda-Guerra identities (discussed in detail in \cite{SKmodel}), as well as the uniform Lipschitz properties in Lemma \ref{VECLemLC1}, one can replace the $\liminf_{N\to\infty} \tilde{A}_{N,M}$ in (\ref{VECFMNANpart22again}) by
$$
f_M^1\bigl(B_\eps(D), \pi\bigr) - f_M^2(\pi)
$$
for some discrete path $\pi\in \Pi_{D_\eps}$ and the functionals defined in (\ref{VECfN1}) and (\ref{VECfN2}) with $\Delta=D_\eps$ due to the constraint in (\ref{VECconstraintDeps}). The discretization introduces some error but it can be made as small as we wish and can be, for example, absorbed in the term $L\eps^{1/4}$ in (\ref{VECFMNANpart22again}). To summarize, the argument in \cite{PPSG} based on the Ghirlanda-Guerra identities and the synchronization property shows that one can find a discrete path $\pi\in \Pi_{D_\eps}$  such that
\begin{align}
\liminf_{N\to\infty} F_N\bigl(\Sigma^N_\eps(D)\bigr)
\geq
\liminf_{M\to\infty} \Bigl( f_M^1\bigl(B_\eps(D), \pi\bigr) - f_M^2(\pi)\Bigr) - L\eps^{1/4}.
\label{VECFMNANpart22again2}
\end{align}
Lemma \ref{VECLem2} then shows that 
$$
\liminf_{M\to\infty} f_M^1\bigl(B_\eps(D),\pi\bigr)
\geq
\inf_{\lambda}\Bigl(  -\sum_{k\leq k'} \lambda_{k,k'} D_{k,k'} + \Phi(\lambda, D_\eps, \pi) \Bigr),
$$
where $\Phi(\lambda, D_\eps, \pi)$ is defined in (\ref{VECPhilambda}) (with $\Delta$ instead of $D_\eps$) and, by (\ref{VECfunp2RPC}),
$$
f_M^2(\pi)
=
\frac{1}{2}
\sum_{0\leq j\leq r-1} x_j  \su\bigl( \theta(\gamma_{j+1}) - \theta(\gamma_{j})\bigr).
$$
Therefore, $\liminf_{N\to\infty} F_N(\Sigma^N_\eps(D))$ is bounded from below by
$$
\inf_{\lambda}\Bigl(  -\sum_{k\leq k'} \lambda_{k,k'} D_{k,k'} + \Phi(\lambda, D_\eps, \pi) \Bigr)
-
\frac{1}{2}
\sum_{0\leq j\leq r-1} x_j  \su\bigl( \theta(\gamma_{j+1}) - \theta(\gamma_{j})\bigr)
-
L\eps^{1/4},
$$
and the proof of Theorem \ref{VECThFELOW} would be finished if we can replace the final constraint $\pi(1) = D_\eps$ in the discrete path $\pi$ by $D$. If we recall the definition of $D_\eps$ in (\ref{VECQLQeps}), it is clear that $D_\eps\leq D$ and $\|D-D_\eps\|_\infty \leq L\sqrt{\eps}.$ Therefore, if we simply extend the path $\pi$ by adding $x_{r+1}=1$ and $\gamma_{r+1}=D$ to the sequences (\ref{VECxs}) and (\ref{VECgammas}), this will modify the above quantity by at most $L\sqrt{\eps}.$ Taking infimum over discrete $\pi\in\Pi_D$ and letting $\eps\downarrow 0$ finishes the proof.


\begin{thebibliography}{999}

\bibitem{AC} Aizenman, M., Contucci, P.: On the stability of the quenched state in mean-field spin-glass models. J. Statist. Phys. \textbf{92}, no. 5-6, 765--783 (1998)

\bibitem{AS2} Aizenman, M., Sims, R., Starr, S.L.: An extended variational principle for the SK spin-glass model. Phys. Rev. B. \textbf{68}, 214403 (2003)

\bibitem{AA} Arguin, L.-P., Aizenman, M.: On the structure of quasi-stationary competing particles systems. Ann. Probab. \textbf{37}, no. 3, 1080--1113 (2009)

\bibitem{AP} Austin, T., Panchenko, D.: A hierarchical version of the de Finetti and Aldous-Hoover representations. Probab. Theory and Related Fields \textbf{159}, no. 3, 809--823 (2014)

\bibitem{MS} Barra, A., Contucci, P., Mingione, E., Tantari, D.: Multi-species mean-field spin-glasses. Rigorous results. Annales Henri Poincar\'e \textbf{16}, no 3, 691--708 (2015)

\bibitem{Bolthausen} Bolthausen, E., Sznitman, A.-S.: On Ruelle's probability cascades and an abstract cavity method. Comm. Math. Phys. \textbf{197}, no. 2, 247--276 (1998) 

\bibitem{Chen-sphere} Chen, W.-K.: The Aizenman-Sims-Starr scheme and Parisi formula for mixed $p$-spin spherical models. Electron. J. Probab. \textbf{18}, no. 94, 1--14 (2013)

\bibitem{ChenChaos0} Chen, W.-K.: Disorder chaos in the Sherrington-Kirkpatrick model with external field. Ann. Probab. \textbf{41}, no. 5, 3345--3391 (2013)

\bibitem{ChenChaos2} Chen, W.-K., Panchenko, D.:  An approach to chaos in some mixed p-spin models. Probab. Theory Related Fields \textbf{151}, no. 1, 389--404 (2013)

\bibitem{ChenChaos} Chen, W.-K.: Chaos in the mixed even-spin models. Comm. Math. Phys. \textbf{328}, no. 3, 867--901 (2014)

\bibitem{Chen15} Chen, W.-K., Hsieh, H.-W., Hwang, C.-R., Sheu, Y.-C.: Disorder chaos in the spherical mean-field model. J. Stat. Phys. \textbf{160}, no. 2, 417--429 (2015)

\bibitem{SCLT} Chen, W.-K., Dey, P., Panchenko, D.: Fluctuations of the free energy in the mixed $p$-spin models with external field. arXiv:1509.07071 (2015)

\bibitem{Dembo} Dembo, A., Zeitouni, O.: Large Deviations Techniques and Applications. Springer-Verlag, New York, (1998)

\bibitem{FL} Franz, S., Leone, M.: Replica bounds for optimization problems and diluted spin systems. J. Statist. Phys. 111, no. 3-4, 535--564 (2003)

\bibitem{FPV} Franz, S., Parisi, G., Virasoro, M.A.: Free-energy cost for ultrametricity violations in spin glasses. Europhysics Letters \textbf{22}, no. 6, 405--411 (1993)

\bibitem{GS} Ghatak, S.K., Sherrington, D.: Crystal field effects in a general $S$ Ising spin glass. J. Phys. C: Solid State Phys. {\bf 10}, 3149 (1977) 

\bibitem{GuerraGG} Guerra, F.: About the overlap distribution in mean field spin glass models. International Journal of Modern Physics B \textbf{10}, no. 13-14, 1675--1684 (1996)

\bibitem{GG} Ghirlanda, S., Guerra, F.: General properties of overlap probability distributions in disordered spin systems. Towards Parisi ultrametricity.  J. Phys. A  \textbf{31}, no. 46, 9149--9155 (1998) 

\bibitem{Guerra} Guerra, F.: Broken replica symmetry bounds in the mean field spin glass model. Comm. Math. Phys. {\bf 233}, no. 1, 1--12 (2003)

\bibitem{MPV} M\'ezard, M., Parisi, G., Virasoro, M.A.: Spin Glass Theory and Beyond. World Scientific Lecture Notes in Physics, 9. World Scientific Publishing Co., Inc., Teaneck, NJ (1987) 

\bibitem{PT} Panchenko, D., Talagrand, M.: Bounds for diluted mean-fields spin glass models. Probab. Theory Related Fields \textbf{130},  no. 3, 319--336 (2004) 

\bibitem{SKcoupled} Panchenko, D.: A note on the free energy of the coupled system in the Sherrington-Kirkpatrick model. Markov Process. Related Fields \textbf{11}, no. 1, 19--36 (2005)

\bibitem{PGS} Panchenko, D.: Free energy in the generalized Sherrington-Kirkpatrick mean field model. Rev. Math. Phys. \textbf{17}, no. 7, 793--857 (2005) 

\bibitem{PGGmixed} Panchenko, D.: The Ghirlanda-Guerra identities for mixed $p$-spin model. C.R. Acad. Sci. Paris, Ser. I \textbf{348}, 189--192 (2010) 

\bibitem{PGG} Panchenko, D.: A connection between Ghirlanda-Guerra identities and ultrametricity.  Ann. of Probab. \textbf{38}, no. 1, 327--347 (2010)

\bibitem{PGG2} Panchenko, D.: Ghirlanda-Guerra identities and ultrametricity: An elementary proof in the discrete case. C. R. Acad. Sci. Paris, Ser. I \textbf{349}, 813--816 (2011)

\bibitem{ACGG} Panchenko, D.: A unified stability property in spin glasses.  Comm. Math. Phys. \textbf{313}, no. 3, 781--790 (2012) 

\bibitem{Pspins} Panchenko, D.: Spin glass models from the point of view of spin distributions.  Ann. of Probab.  {41}, no. 3A, 1315--1361 (2013)

\bibitem{PUltra} Panchenko, D.: The Parisi ultrametricity conjecture. Ann. of Math. (2) \textbf{177}, no. 1, 383--393 (2013)

\bibitem{SKmodel} Panchenko, D.: The Sherrington-Kirkpatrick Model. Springer Monographs in Mathematics. Springer-Verlag, New York (2013)

\bibitem{PPF} Panchenko, D.: The Parisi formula for mixed $p$-spin models. Annals of Probability \textbf{42}, no. 3, 946--958 (2014)

\bibitem{PMS} Panchenko, D.: The free energy in a multi-species Sherrington-Kirkpatrick model. Annals of Probability \textbf{43}, no. 6, 3494--3513 (2015)

\bibitem{HEPS} Panchenko, D.: Hierarchical exchangeability of pure states in mean field spin glass models. Probab. Theory and Related Fields \textbf{161}, no. 3, 619--650 (2015)

\bibitem{PTChaos} Panchenko, D.: Chaos in temperature in generic $2p$-spin models.  arXiv:1502.03801 (2015)

\bibitem{PPSG} Panchenko, D.: Free energy in the Potts spin glass. arXiv:1512.00370 (2015)

\bibitem{Parisi79} Parisi, G.: Infinite number of order parameters for spin-glasses. Phys. Rev. Lett. \textbf{43}, 1754--1756 (1979)

\bibitem{Parisi} Parisi, G.: A sequence of approximate solutions to the S-K model for spin glasses. J. Phys. A \textbf{13}, L-115 (1980) 

\bibitem{RC} Rizzo, T., Crisanti, A.: Chaos in temperature in the Sherrington-Kirkpatrick model. Phys. Rev. Lett. \textbf{90}, 137201 (2003) 

\bibitem{Ruelle} Ruelle, D.: A mathematical reformulation of Derrida's REM and GREM.  Comm. Math. Phys.  {\bf 108},  no. 2, 225--239 (1987)

\bibitem{SK} Sherrington, D., Kirkpatrick, S.: Solvable model of a spin glass. Phys. Rev. Lett. {\bf 35}, 1792--1796 (1975)

\bibitem{SG} Talagrand, M.:  Spin Glasses: a Challenge for Mathematicians.  Ergebnisse der Mathematik und ihrer Grenzgebiete. 3. Folge A Series of Modern Surveys in Mathematics, Vol. 43. Springer-Verlag (2003) 

\bibitem{TPF} Talagrand, M.: The Parisi formula. Ann. of Math. (2) \textbf{163}, no. 1, 221--263 (2006)

\bibitem{T-sphere} Talagrand, M.: Free energy of the spherical mean field model. Probab. Theory Related Fields  \textbf{134}, no. 3, 339--382 (2006)

\bibitem{PM} Talagrand, M.: Parisi measures.  J. Funct. Anal. \textbf{231}, no. 2, 269--286 (2006)

\bibitem{TalUltra} Talagrand, M.: Mean field models for spin glasses: some obnoxious problems. Lecture Notes in Mathematics, Vol. 1900, 63--80 (2007).
 
\bibitem{SG2} Talagrand, M.: Mean-Field Models for Spin Glasses. Ergebnisse der Mathematik und ihrer Grenzgebiete. 3. Folge A Series of Modern Surveys in Mathematics, Vol. 54, 55. Springer-Verlag (2011)

\end{thebibliography}
\end{document}